\documentclass[a4paper,reqno]{amsart}
\usepackage{mathrsfs,amssymb,stmaryrd,bbm,mathabx}

\usepackage{amsmath}
\usepackage{mathtools}
\usepackage
{hyperref}
\usepackage[textsize=scriptsize,color=orange]{todonotes}
\usepackage[inline,shortlabels]{enumitem}
\usepackage{datetime2}
\usepackage{graphicx}
\usepackage{tikz-cd}
\usepackage{changebar}
\usepackage{leftidx}
\DeclareMathAlphabet{\mathpzc}{OT1}{pzc}{m}{n}

\theoremstyle{plain}
\newtheorem{thm}{Theorem}
\newtheorem{cor}[thm]{Corollary}
\newtheorem{lemma}[thm]{Lemma}

\theoremstyle{remark}

\newtheorem{rk}[thm]{Remark}

\newtheorem{exa}[thm]{Example}

\theoremstyle{definition}

\newcommand{\ti}{\mbox{--}}
\newcommand{\Ob}{\mathrm{Ob}}
\newcommand{\A}{\ensuremath{\mathcal{A}}}
\newcommand{\B}{\ensuremath{\mathcal{B}}}
\newcommand{\C}{\ensuremath{\mathcal{C}}}
\newcommand{\D}{\ensuremath{\mathcal{D}}}
\newcommand{\E}{\ensuremath{\mathcal{E}}}
\newcommand{\G}{\ensuremath{\mathcal{G}}}
\newcommand{\Z}{\ensuremath{\mathcal{Z}}}
\newcommand{\M}{\ensuremath{\mathcal{M}}}

\newcommand{\Mod}{\mathsf{Mod}}
\renewcommand{\mod}{\mathsf{mod}}

\newcommand{\Vect}{\ensuremath{\mathbf{Vect}}}
\newcommand{\vect}{\ensuremath{\mathbf{vect}}}

\newcommand{\xtor}[1]{\cdl[@1]{{} \ar[r]|-{\object@{|}}^{#1} & {}}}

\newcommand{\adj}{\dashv}
\newcommand{\opp}{\mathrm{op}}
\newcommand{\Set}{\ensuremath\mathbf{Set}}
\newcommand{\Hom}{\mathrm{Hom}}
\newcommand{\End}{\mathrm{End}}
\newcommand{\Fun}{\mathrm{Fun}}
\newcommand{\Lan}{\mathrm{Lan}}
\newcommand{\id}{\mathrm{id}}
\newcommand{\eps}{\varepsilon}
\newcommand{\Nat}{\mathrm{Nat}}
\newcommand{\Aut}{\mathrm{Aut}}
\newcommand{\rhom}[2]{[#1,#2]^r}
\newcommand{\Lhom}[2]{[#1,#2]^l}
\newcommand{\lhom}[2]{[#1,#2]}

\newcommand{\pLhom}[2]{\underline{[#1,#2]}^l}
\newcommand{\plhom}[2]{\underline{[#1,#2]}}
\newcommand{\iso}{\xrightarrow{\cong}}

\newcommand{\cp}{\rtimes}

\newcommand{\dar}[2]{\ar@<2pt>[r]^-{#1}\ar@<-2pt>[r]_-{#2}}
\newcommand{\blabla}[1]{\quad \text{#1}\quad}
\newcommand{\kk}{\Bbbk}
\newcommand{\SL}{\mathrm{sl}}
\newcommand{\ALPHA}{\boldsymbol{\alpha}}
\newcommand{\BETA}{\boldsymbol{\beta}}
\newcommand{\ZETA}{\boldsymbol{\zeta}}
\newcommand{\THETA}{\boldsymbol{\theta}}
\newcommand{\XI}{\boldsymbol{\xi}}
\newcommand{\PSI}{\boldsymbol{\Psi}}

\newcommand{\un}{\mathbbm{1}}

\newcommand{\pt}{{\bf *}}

\newcommand{\prettyadj}[4]
{\xymatrix@1{
    #3
    \ar@<4pt>[r]^-{#1}
    \ar@<-4pt>@{<-}[r]_-{#2}
    \ar@{}[r]|-\bot
    &
    #4}}
\newcommand{\ldual}[1]{\leftidx{^\vee}{\!#1}{}}

\newcommand{\alf}{\mathpzc{a}}
\newcommand{\bee}{\mathpzc{b}}
\newcommand{\Es}{\mathpzc{S}}
\newcommand{\inv}{\times}

\DeclareMathOperator*{\sstimes}{\scalebox{.9}{$\times$}}
\DeclareMathOperator*{\sscoprod}{\scalebox{.9}{$\coprod$}}

\begin{document}
\title{Slack Hopf monads}
\author[A.~Brugui\`eres]{Alain Brugui\`eres}
\address{Institut Montpelliérain Alexander Grothendieck, Universit\'e de Montpellier, CNRS, CC 051 - Place Eug{\`e}ne Bataillon, 34095 Montpellier CEDEX 5, France}
\email{alain.bruguieres@umontpellier.fr}
\author[M.~Haim]{Mariana Haim}
\address{Centro de Matem\'atica, Facultad de Ciencias, Universidad de la
  Rep\'ublica. Igu\'a 4225, Montevideo, Uruguay.}
\email{negra@cmat.edu.uy}
\author[I.~L\'opez~Franco]{Ignacio L\'opez Franco}
\address{Departamento de Matem\'atica y Aplicaciones, CURE, Universidad de la
  Rep\'ublica. Tacuaremb\'o s/n, Maldonado, Uruguay.}
\email{ilopez@cure.edu.uy}
\thanks{The authors acknowledge the support of Universidad de la Rep\'ublica, PEDECIBA and the IFUMI (Franco-Uruguayan Institute for Mathematics and  Interactions).}


\begin{abstract}
  Hopf monads generalise Hopf algebras. They clarify several aspects of the theory of Hopf algebras and capture several related structures such as weak Hopf algebras and Hopf algebroids. However,
  important parts of Hopf algebra theory are not reached by Hopf
  monads, most noticeably  Drinfeld's quasi-Hopf algebras. In this paper
  we introduce a generalisation of Hopf monads, that we call slack Hopf
  monads. This generalisation retains a clean theory and is flexible enough to
  encompass quasi-Hopf algebras as examples. A slack Hopf monad is a colax magma monad $T$ on a magma category $\C$ such that the forgetful functor $U^T \colon \C^T \to \C$ `slackly' preserves internal Homs.
  We give a number of different descriptions of slack Hopf monads, and study special cases such as slack Hopf monads on cartesian categories and $\kk$-linear exact slack Hopf monads on $\Vect_\kk$, that is comagma algebras such that a modified fusion operator is invertible. In particular, we characterise quasi-Hopf algebras in terms of slackness.
\end{abstract}
\maketitle
\setcounter{tocdepth}{1}
\tableofcontents

\section*{Introduction}
\label{sec:introduction-1}
Hopf monads were introduced at first in the context of autonomous categories \cite{BruVir}, with a view on quantum topological applications, and later extended to arbitrary monoidal categories
in \cite{MR2793022}.
A comonoidal monad (called Hopf monads in \cite{Moer}, opmonoidal monads in
\cite{zbMATH01863401} and bimonads in \cite{BruVir,MR2793022})
on a monoidal category $\C$ is a monoid in the category of comonoidal
endofunctors of $\C$. If $T$ is a comonoidal monad on $\C$, the Eilenberg-Moore
category $\C^T$ of $T$-modules in $\C$ is monoidal, and the forgetful functor
$U^T\colon \C^T \to \C$ is monoidal strict \cite{Moer}.
Comonoidal monads generalise bialgebras, and by the same token, Hopf monads generalise Hopf algebras. There are however several characterisations of Hopf algebras, which lead to different definitions for Hopf monads.  Thus in \cite{BruVir}, a left Hopf monad on a left autonomous category $\C$ is defined as a comonoidal monad $T$ such that $\C^T$ is left autonomous, a condition equivalent to the existence of a certain natural transformation $S^l_X : T(\ldual{T(X)}) \to \ldual{X}$ called a left antipode; whereas in \cite{MR2793022}, a left Hopf monad on an arbitrary monoidal category $\C$ is a comonoidal monad $T$ such that the left fusion operator
$$H^l_{X,Y} = (TX \otimes \mu_Y)T_2(X,TY) : T(X \otimes TY) \to TX \otimes TY$$
is an isomorphism (where $\mu$ and $T_2$ denote respectively the multiplication and the comonoidal structure of $T$). This definition coincides with the previous one when $\C$ is left autonomous and, if $\C$ is left monoidal closed, means that $\C^T$ is left monoidal closed and $U^T$ preserves left internal Homs.\\
In this full-fledged version, Hopf monads capture not only Hopf algebras, but many variants thereof, such as weak Hopf algebras, Hopf algebroids, Hopf algebras in braided categories, and so on.

However, the theory of Hopf monads fails to encompass some important aspects of Hopf algebra theory, and such is the case, most notably, of Drinfeld's quasi-Hopf algebras --- for an obvious reason: if $A$ is a quasi-bialgebra over a field $\kk$, its category of modules $A\ti\Mod$ is monoidal but the forgetful functor $A\ti\Mod \to \Vect_\kk$ is not monoidal, unlike the forgetful functor of a comonoidal monad.
This difficulty can be circumvented by modifying the notion of a comonoidal monad so as to allow for a `quasi-monoidal' forgetful functor by introducing an associator.

There is however a more fundamental difficulty when it comes to generalising quasi-Hopf algebras: the `quasi-Hopf' property (the existence of a quasi-antipode $(S,\alf,\bee)$) is a purely algebraic condition and does not readily lend itself to a categorical characterisation. In fact, such a characterisation exists in the case of profinite quasi-bialgebras: a profinite quasi-bialgebra is quasi-Hopf if and only if its category of finite dimensional modules is left autonomous, and the forgetful functor is left autonomous, that is, preserves left duals (see Theorem~\ref{thm:quasi-auton});
but this does not seem to hold for arbitrary quasi-Hopf algebras. Therefore, in order to exploit this partial characterisation of quasi-Hopf, we would need to restrict ourselves to monoidal categories having `enough duals' and monads satisfying certain technical conditions. This is not a path we are willing to tread in this article, as we would have to discard for instance cartesian categories, where only the unit object has duals.

We can also not bypass assumptions on the existence of duals by requiring the invertibility of the fusion operator, because the fusion operator $(A \otimes m)(\Delta \otimes A): A \otimes A \to A \otimes A$ of a quasi-Hopf algebra $A$ (with product $m$ and coproduct $\Delta$) is not necessarily bijective. It has been noted by Drinfeld \cite{QHA} and Schauenburg \cite{Schau:CFQHA} that for a quasi-Hopf algebra, there is a sort of modified fusion operator which is bijective, but this modified fusion operator is constructed from the fusion operator using a quasi-antipode and therefore cannot serve as a definition of `quasi-Hopf'. On the flip side, Schauenburg has noted that the existence of this modified fusion operator is closely related to an interesting fact: if $A$ is a quasi-Hopf algebra, then $A\ti \Mod$ is monoidal left closed, and the forgetful functor preserves the left internal Hom functors $\lhom{V}{W}$, but not necessarily the evaluation morphisms $\lhom{V}{W} \otimes V \to W$ (see Section~\ref{sec:quasi-hopf-slack}).

So, what does Hopf mean? It seems that we can't answer this question without making choices, and in this paper we decided to explore a very relaxed interpretation, namely in terms of slack preservation of internal Homs. This lead us to a radical choice: we work in the world of magma categories. A magma category is a category $\C$ equipped with an arbitrary functor $\otimes : \C \times \C \to \C$.  In this desolate landscape, a few denizens of the lush monoidal garden survive. Colax magma functors are a substitute for comonoidal functors, and colax adjunctions behave pretty much like comonoidal adjunctions.  A colax magma monad $T$ on a magma category $\C$ is characterised by the fact that $\C^T$ is a magma category and $U^T$ is a strict magma functor. Duals do not belong here, but we can still define left and right internal Homs.

What does it mean for a functor $U: \D \to \C$ between magma categories to preserve internal Homs? Assume for simplicity that $\C$ and $\D$ have left internal Homs.\footnote{This assumption can be lifted at the price of a reformulation in terms of presheaves.} If $U$ is a colax magma functor, there is for $A, B$ in $\D$ a canonical morphism $U(\lhom{A}{B}) \to \lhom{U(A)}{U(B)}$
and we may request this morphism to be an isomorphism. However, this form of preservation of internal Homs, the one satisfied by the forgetful functor of a Hopf monad, would be too strong for our purposes, as it is not satisfied by the forgetful functor $A\ti \Mod \to \Vect_\kk$ of a general quasi-bialgebra $A$. Instead, we say that a functor $U : \D \to \C$ slackly preserves left internal Homs, or in short $U$ is slack left closed, if it is equipped with an arbitrary natural isomorphism
$U(\lhom{A}{B}) \iso \lhom{U(A)}{U(B)}$. Such an isomorphism is called a slack left closed structure on $U$.
We define an adjunction $F \adj U : \C \to \D$ to be slack left Hopf if $U: \D \to \C$ is slack left closed. There is a bijection between slack left closed structures on $U$ and slack left Frobenius isomorphisms $\Psi$ on $F \adj U$, that is, natural isomorphisms $\Psi_{X,A} : F(X \otimes U(A)) \iso F(X) \otimes A$ for $X \in \C$, $A \in \D$ (see Section~\ref{sec-slack-Hopf-adj}).

A slack left Hopf monad on a magma category $\C$ is a colax magma monad $T$ whose forgetful functor $U^T$ is slack left closed. This is equivalent to the existence of a natural transformation
$$\beta_{X,Y} : X \otimes Y \to TX \otimes TY$$
such that the natural morphism $$\Phi^\beta_{X,Y} = (\mu_X \otimes \mu_Y)T_2(TX,TY)T((TX \otimes \mu_Y)\beta_{X,TY}) : T(X \otimes TY) \to TX \otimes TY$$
is an isomorphism (Theorem~\ref{thm-ll-psi-beta-2}). Such a $\beta$ is called a slack (left) Hopf structure on $T$.
In general there are many slack Hopf structures, but they form a unique orbit under the free action of a certain group (see Section~\ref{sec:unicity}). If $\C$ is closed, a slack Hopf structure allows one to lift the internal Homs of $\C$ to $\C^T$, which is therefore also closed, and the lifted internal Homs can be defined in terms of a natural transformation which plays the r\^ole of the (left) antipode of \cite{MR2793022}.
In other words the forgetful functor of a slack Hopf monad not only (slackly) preserves internal Homs, but it creates them as well (Theorems~\ref{thm-ll-mon-mod} and~\ref{thm-antipode}).

A Hopf monad on a monoidal category $\C$ is a special example of a slack Hopf monad, with trivial slack Hopf structure $\beta = \eta \otimes \eta$, $\eta$ being the unit of $T$; for this choice of $\beta$, $\Phi^\beta$ is the left fusion operator $H^l$. The notion of slack Hopf is so weak that we expect comonoidal monads can be slack Hopf without being Hopf. However, we do not have so far any example of such a behaviour, and in the case of a cartesian category, we show that any slack Hopf comonoidal monad is a Hopf monad (Theorem~\ref{th-hopf-cart}). In particular a small category which is slack Hopf, viewed as a colax monad, is a groupoid (Theorem~\ref{thm-cat-monad}), and it is also the case, for different reasons, when it is viewed as a lax comonad (Theorem~\ref{thm-cat-comonad}).

Comagma algebras are an important source of examples. A comagma algebra over $\kk$ is a  $\kk$-algebra $A$ equipped with an algebra morphism $\Delta : A \to A \otimes A$. The monad $A \otimes ?$ on $\Vect_\kk$ is a colax magma monad, and indeed any $\kk$-linear exact colax magma monad on $\Vect_\kk$ is of this type. A slack left Hopf structure on such a monad corresponds with an element $v \in A \otimes A$ such that the map:
$$H^v: A \otimes A \to A \otimes A, \quad x \otimes y \mapsto x_{(1)}v^{(1)} \otimes x_{(2)}v^{(2)}y$$
is a bijection. Such a $v$ is also called a slack Hopf structure.

If $A$ is a bialgebra, we characterise the slack left Hopf structures $v$ such that $A$ is Hopf in the profinite case (Theorem~\ref{thm-slack-Hopf-bialg}). A commutative bialgebra which is slack left Hopf is a Hopf algebra (Corollary~\ref{cor:commut bialg}).

If $A$ is a quasi-Hopf algebra, a quasi-antipode $(\Es,\alf,\bee)$ gives rise to a slack left Hopf structure $v =\phi^{(-1)}\otimes  \phi^{(-2)}\,\bee \, \Es( \phi^{(-3)})$, and the fact that $H^v$ is an isomorphism was already noted by Schauenburg. We call such a $v$ associated to a quasi-antipode a left Hopf structure.
If $A$ is a (non-necessarily profinite) quasi-bialgebra, we define the slackness $\SL(v)$ of a slack left Hopf structure $v$, which is an element of the enveloping algebra of $A$ measuring the obstruction to $v$  being a Hopf structure, and so characterise the left Hopf structures among slack left Hopf structures (Theorem~\ref{thm-crit-quanti}), that, is we give a criterion for a slack left Hopf quasi-bialgebra to be quasi-Hopf (Corollary~\ref{cor-quanti}).
This allows us to give a purely categorical characterization of a quasi-Hopf algebra: it is a quasi-bialgebra
whose forgetful functor admits a slack left closed structure which preserves comparison morphisms (Theorem~\ref{thm-char-quasi}).
In a quasi-Hopf algebra $A$  we show that the set of slack left Hopf structures forms a torsor under the group of invertibles of the enveloping algebra of $A$, inside of which left Hopf structures form an orbit under $A^\inv$ (Theorem~\ref{thm:quanti}).

In spite of the fact that we dispose of these criteria, we still do not have examples of slack Hopf bialgebras or quasi-bialgebras which are not Hopf. On the other hand, any monoid in a monoidal category defines a slack Hopf monad (Example~\ref{ex:alg-slack}) and in particular any $\kk$-algebra $A$ has a structure of a slack Hopf comagma algebra, and for $n \ge 2$, $M_n(\kk)$ has a structure of a slack Hopf comagma algebra for which the fusion operator $(A \otimes m)(\Delta\otimes A)$ is not invertible --- but $M_n(\kk)$ is not a bialgebra (Example~\ref{exa:mat}).

\subsection*{Organisation of the article} In Section~\ref{sec:conventions} we recall certain classical facts about categories, Kan extensions, coends, presheaves and monads, and fix some notation.
Section~\ref{sec:main-results} deals with magma categories, colax magma functors, colax magmatic adjunctions, colax magma monads, internal Homs and magma closed categories. The results of this section closely mimick analogous classical results for the monoidal counterparts of the objects considered, their sole interest residing in the constatation that they somehow manage to survive in the inhospitable world of magma categories.
In Section~\ref{sec:slack-closed}, we define what it means for a functor to slackly preserve internal Homs, which leads to the notions of slack closed functor and slack Hopf adjunction.
Section~\ref{sec:slack-hopf-adjunct} is devoted to slack Hopf monads and their basic properties, slack Hopf structures, slack antipodes and the construction of internal Homs in the category of modules of a slack Hopf monad. In Section~\ref{sec-unit-cart} we consider the case where the magma category $\C$ has a (one-sided) unit object and we show two technical results about slack Hopf structures which serve to prove the main theorem of the section, namely that on a cartesian category, slack Hopf monoidal monads coincide with Hopf monads. In 
particular we tackle the example of a monad or a comonad associated with a small category.
Section~\ref{sec:examples} is devoted to comagma monoids and comagma algebras, and Section~\ref{sec:qha} to the special case of quasi-bialgebras and quasi-Hopf algebras. In a short conclusion (Section~\ref{sec:concl}) we address two questions which remain open.

\section{Preliminaries and notations}
\label{sec:conventions}
All categories are assumed to be locally small, unless otherwise specified. Given two categories $\A$, $\B$, we denote by $[\A,\B]$ the category of functors $\A \to \B$ and natural transformations (which may not be locally small).

For clarity, when we work with categories $\C$ and $\D$, as it will often be the case, we will usually denote objects of $\C$ by $X,Y,...$ and objects of $\D$ by $A, B, ...$. This convention will hopefully help decypher complicated formulae.

 Given two adjunctions $G\adj V$, $G' \adj V'$ and two functors  $M$, $N$
forming a diagram:
$$
\begin{tikzcd}
    \C'\ar[d,"M"']
    \ar[r,shift left=1.2ex,"{G'}"]
    \ar[r,phantom,"\bot"]
    &
    \D'\ar[l,shift left=1ex,"{V'}"]
    \ar[d,"N"]
    \\
    \C\ar[r,shift left=1.5ex,"G"]
    \ar[r,phantom,"\bot"]
    &
    \D\ar[l,shift left=1ex,"V"]
\end{tikzcd}
$$
we have a natural bijection
$\Nat(GM,NG') \iso \Nat(MV',VN)$. A natural transformation $\alpha$ on the left hand side and the corresponding natural transformation $\beta$ on the right hand side are said to be each other's \emph{mate}. Denoting by $h: \id \to VG$, $e: GV \to \id$, $h': \id \to V'G'$, $e': G'V' \to \id$ the adjunction morphisms, the mates $\alpha$ and $\beta$ are related by:
$$\beta =(VNe')(V\alpha V')(hMV') \blabla{and} \alpha = (e NG')(G \beta G')(GMh').$$

\subsection{Left Kan extensions, coends.}
Given two functors $U \colon \D \to \C$, $F \colon \D \to \E$ between locally small categories, a \emph{left Kan extension of $F$ along $U$} is a functor $H : \C \to \E$ equipped with a natural transformation $\eta \colon F \to HU$
satisfying the following universal property: for every functor $H': \C \to \E$ and every natural transformation $\alpha \colon F \to H'U$, there exists a unique natural transformation $\widetilde{\alpha} \colon H \to H'$ such that $\alpha_A = \widetilde{\alpha}_{UA}\,\eta_A$ for $A$ in $\D$.

A left Kan extension of $F$ along $U$, if it exists, is unique up to unique isomorphism, and is denoted by $\Lan_U (F)$. (See \cite{CWM} for properties of Kan extensions.)

If $\E = \Set$, the left Kan extension $\Lan_U (F)$ can be interpreted in terms of coends. More precisely, $\Lan_U(F)$ exists if and only if the coend $\int^{A \in \D} \C(UA,X) \sstimes F(A)$ exists for every $X$ in $\C$, and in that case it is isomorphic to the aforementioned coend.

\subsection{Cocontinuous extensions.} Given a locally small category $\C$, we denote by $\widehat{\C}$ the category $\Fun(\C^\opp,\Set)$ of presheaves on $\C$ (which may not be locally small). Recall that $\C$ can be viewed as a full subcategory of $\widehat{\C}$ via the Yoneda embedding $X \mapsto \C(-,X)$.
A functor $U : \D \to \C$ between locally small categories defines by precomposition a functor $U^* : \widehat{\C} \to \widehat{\D}$.

A \emph{cocontinuous extension of $U$} is a functor $U_! \colon \widehat{\D} \to \widehat{\C}$ which coincides with $U$ on $\D \subset \widehat{\D}$ and preserves small limits. If it exists, it is unique up to unique isomorphism.
The following conditions are equivalent:
\begin{enumerate}[(i)]
\item $U$ admits a cocontinuous extension $U_! : \widehat{\D} \to \widehat{\C}$;
\item $U^*$ admits a left adjoint;
\item the coend $\int^{A \in \D} \C(-,U(A))\sstimes X(A)$ exists for every presheaf $X$ on $\D$.
\end{enumerate}
If these conditions hold then the cocontinuous extension $U_!$ is left adjoint to $U^*$, and it is defined by the coend
$$U_!(X) =  \int^{A \in \D} \C(-,U(A))\sstimes X(A) = \Lan_{U^\opp} (X).$$

Note that these conditions hold if for each object $c$ of $\C$, the functor $$\C(c,U(-)) : \D \to \Set$$ is a small colimit of representables. In particular, they hold
\begin{itemize} \item if $U$ has a left adjoint $F$, and in that case $U_! = F^*$;
\item if $\D$ is essentially small.
\end{itemize}

\subsection{Monads}
Let  $\C$ be a category. A \emph{monad on $\C$} is a monoid in the monoidal categories of endofunctors of $\C$, that is, a triple $(T,\mu,\eta)$ where $T$ is an endofunctor of $\C$, $\mu: T^2 \to T$ and $\eta: \id_\C \to T$ are natural transformations satisfying the following associativity and unity axioms:
$$\mu(\mu \sstimes \C) = \mu(\C \sstimes \mu), \qquad \mu(\eta \sstimes \C) = \id_\C = \mu(\C \sstimes \eta).$$
Given such a monad $T$, a \emph{$T$-module \footnote{In the literature $T$-modules in $\C$ are often called $T$-algebras.} in $\C$} is an object $X$ of $\C$ endowed with an \emph{action of $T$}, that is a morphism $r : TX \to X$ satisfying the following conditions:
$$r \mu_X = r T(r),  \qquad r\eta_X = \id_X.$$
Given two $T$-modules $(X,r)$, $(Y,\rho)$ in $\C$, a $T$-module morphism $f : (X,r)\to (Y,\rho)$ is a morphism $f : X \to Y$ in $\C$ satisfying the $T$-linearity condition: $f r= \rho T(f)$.

In this way $T$-modules in $\C$ form a category $\C^T$, called the Eilenberg-Moore category of $T$. The  forgetful functor
$U^T : \C^T \to \C$ admits a left adjoint, the free module functor $F^T : \C \to \C^T$, defined on objects, by $F^T(X) = (TX,\mu_X)$ and on morphisms, by $F^T(f) = T(f)$, so that $T = U^TF^T$.

Any adjunction $F \adj U : \D \to \C$, with unit $\eta$ and counit $\varepsilon$, gives rise to a monad $T = UF$ on $\C$, with multiplication $\mu = U\varepsilon_U$
and unit $\eta$.

\subsection{Convolution product} \label{sec:convo}
Given a monad $T$ on a category $\C$, the adjunction bijection $$\Nat(\id_\C, T) \iso \End(U^T), \quad a \mapsto \check a = U^T(\eps)a_{U^T}$$
is an isomorphism of monoids when $\Nat(\id_C, T)$ is equipped with the \emph{convolution product}: for $a,b  \in \Nat(\id_\C, T)$,
$$a * b = \mu T(b) a = \mu a_T b,$$
whose unit is $\eta$.
The other adjunction bijection $$\Nat(\id_\C,T)\iso\End(F^T), \quad a \mapsto \hat a = \eps_{F^T} {F^T}(a),$$ is an isomorphism when $\Nat(\id_\C,T)$ is equipped with the \emph{Kleisli composition} \cite{zbMATH03223713}, that is, the opposite of the convolution product. 

\subsection{Lifting data.}
Let $\C$, $\D$ be categories, $(T,\mu,\eta)$ and $(P,m,u)$ be monads on $\C$ and $\D$ respectively, and let $G: \C \to \D$ be a functor.
Then there is a one-to-one correspondance between
\begin{enumerate}
\item functors $\widetilde{G} : \C^T \to \D^P$ such that $U^P\widetilde{G} = G U^T$ (called \emph{liftings of $G$})
\item natural transformations : $\zeta : P\circ G \to G \circ T$ satisfying:\begin{equation}\label{eq:0-1}
\zeta m_G = G(\mu) \zeta_T P(\zeta), \qquad \zeta u_G = G(\eta).
\end{equation}
\end{enumerate}
The correspondence between $\widetilde{G}$ and $\zeta$ is as follows:
$$\widetilde{G}(X,r) = (GX,G(r)\zeta_X) \blabla{and} \zeta_X = U^P\eps_{\widetilde{G}F^TX}PG(\eta_X).$$
A natural transformation $\zeta$ satisfying (\ref{eq:0-1}) is called a \emph{lifting datum for $G$} and the corresponding lifting of $G$ is denoted by $\widetilde{G}^\zeta$.

Given two functors $G , H : \C \to \D$ and two lifting data $\zeta : P\circ G \to G \circ T$, $\xi: P\circ H \to H \circ T$, there is a one-to-one correspondance between:
\begin{enumerate}
\item natural transformations $\widetilde{\gamma} : \widetilde{G}^\zeta \to \widetilde{H}^\xi$;
\item natural transformations $\gamma: G \to H\circ T$ satisfying:
\begin{equation}
H(\mu)\gamma_T \zeta = H(\mu) \xi_T  P(\gamma),
\end{equation}
\end{enumerate}
the relation between $\widetilde\gamma$ and the corresponding $\gamma$ being:
$$\widetilde{\gamma}_{(X,r)} = H(r)\gamma_X \blabla{and} \gamma_X = \widetilde\gamma_{F^TX} G(\eta_X).$$

\section{Magma categories}
\label{sec:main-results}

\subsection{Magma categories}
\label{sec:lax magma-categories}

A \emph{magma category} is a magma in the category of categories, that is, a category $\C$ endowed with a functor $\otimes : \C \sstimes \C \to \C$.

Given two magma categories $\C$, $\D$, a \emph{lax magma functor} $F : \C \to \D$ consists in a functor $F : \C \to \D$ and a natural transformation
$$F^2(X,Y): F(X) \otimes F(Y) \to F(X\otimes Y), \quad X, Y \in \C.$$

A lax magma functor $F$ is said to be \emph{strong} (resp. \emph{strict}) if $F^2$ is an isomorphism (resp. an identity).

A natural transformation $\alpha : F \to G$ between lax magma functors is \emph{lax magmatic} if it satisfies:
$$\alpha_{X\otimes Y}F^2(X,Y)= G^2(X,Y)(\alpha_{X} \otimes \alpha_{Y}).$$
This defines a 2-category, whose $0$-cells are magma categories, $1$-cells are lax magma functors and $2$-cells are lax magmatic natural transformations.

Dually, a \emph{colax magma functor} $F : \C \to \D$ is a functor $F : \C \to \D$ equipped with a natural transformation
$$F_2(X,Y) : F(X \otimes Y) \to F(X) \otimes F(Y), \quad X, Y \in \C.$$
A natural transformation $\alpha : F \to G$ between colax magma functors is \emph{colax magmatic} if it satisfies:
$$G_2(X,Y)\alpha_{X \otimes Y} = (\alpha_X \otimes \alpha_Y) F_2(X,Y).$$
This defines a 2-category, whose $0$-cells are magma categories, $1$-cells are colax magma functors and $2$-cells are colax magmatic natural transformations.

A \emph{strong magma functor} is a colax magma functor $F$ such that $F_2$ is an isomorphism\footnote{Or equivalently a lax magma functor $F$ such that $F^2$ is an isomorphism, but we will deal mostly with colax magma functors.}

\subsection{Colax magmatic adjunctions and colax magma monads.}

A \emph{colax magmatic adjunction} is an adjunction $F \adj U$, where $F : \C \to \D$ and $U : \D \to \C$ are colax magma functors between magma categories, and the unit $\eta: \id_\C \to UF $ and counit $\eps: FU \to \id_\D$ of the adjunction are colax magmatic natural transformations.

\begin{lemma} Consider an adjunction $F \adj U$ between magma categories.
\begin{enumerate}
\item There is a one-to-one correspondance between colax magma structures $F_2$ on $F$ and lax magma structures $U^2$ on $U$, given by:
$$F_2(X,Y)= \eps_{FX \otimes FY}F(U^2(FX,FY)) F(\eta_X \otimes \eta_Y)$$
$$U^2(A,B) = U(\eps_A \otimes \eps_B)U(F_2(UA,UB))\eta_{UA \otimes UB}$$
\item If the adjunction $F\adj U$ is colax magmatic, then $U$ is a strong magma functor;
\item If $U$ is a strong magma functor, then $F$ admits a unique colax magma structure making the adjunction $F\adj U$ colax magmatic.
\end{enumerate}
\end{lemma}

\begin{proof} This can be viewed as a consequence of Kelly's theorem on doctrinal adjunctions (\cite{kelly_doctrinal_1974}).
The bijection of Assertion (1) is the adjunction bijection:
$$\Nat(F \circ \otimes_\C, \otimes_\D \circ (F \sstimes F)) \cong \Nat(\otimes_\C \circ (U \sstimes U), U \circ \otimes_\D),$$
so that a colax magmatic structure $F_2$ on $F$ and the corresponding lax magmatic structure  $U^2$ on $U$ are related by:
 \begin{align*}
 UF_2(X,Y)\eta_{X\otimes Y} &= U^2(FX,FY)(\eta_X \otimes \eta_Y),\\
\eps_{A \otimes B}FU^2(A,B) &=(\eps_A \otimes \eps_B)F_2(UA \otimes UB).
\end{align*}

One deduces from these identities firstly that if $F \adj U$ is colax magmatic the lax magma structure $U^2$ corresponding with $F_2$ is inverse to $U_2$, that is, Assertion~(2), and secondly, that
if $U$ is a strong magma functor, a colax magma structure $F^2$ on $F$ makes $F\adj U$ a colax magmatic adjunction if and if only if it is the one associated with $U^2 = U^{-1}_2$, that is, Assertion~(3).
\end{proof}

\subsection{Magma equivalences}

A \emph{magma equivalence} $F : \C \to \D$ between two magma categories $\C$, $\D$ is a colax magma functor $F : \C \to \D$ whose underlying functor is an equivalence of categories. It follows from the previous lemma that $F$ is then automatically a strong magma functor.

\subsection{Colax magma monads.}

A \emph{colax magma monad} on a magma category $\C$ is a monad $(T,\mu,\eta)$ on $\C$ endowed with a colax magma structure $T_2$ on $T$
such that $\mu : T^2 \to T$ and $\eta: \id_\C \to T$ are colax magmatic.

\begin{exa} If $\C$, $\D$ are magma categories and $U : \C \to \D$ is a strong magma functor having a left adjoint $F$, then $F \adj U$ is a colax magmatic adjunction and its monad $T = UF$ is a colax magma monad.
\end{exa}

In fact all colax magma monads are of this form. More precisely, if $T$ is a colax magma monad on a magma category $\C$, its Eilenberg-Moore category $\C^T$ is a magma category and the forgetful functor $U^T : \C^T \to \C$ is a strict magma functor. Its left adjoint $F^T : \C \to \C^T$ is therefore a colax magma functor, and the adjunction $F^T \adj U^T$ is colax magmatic.

\begin{exa}
A comonoidal monad on a monoidal category is a colax magma monad.
\end{exa}

\begin{thm} Let $\C$ be a magma category and $T$ be a monad on $T$. There exists a bijection between the following:
\begin{enumerate}[(i)]
\item magma structures on $\C^T$ such that $U^T$ is strict colax magmatic;
\item colax magma monad structures on $T$.
\end{enumerate}
\end{thm}

\begin{proof}
A magma structure on $\C^T$ such that $U^T$ is strict colax magmatic is a functor $\widetilde{\otimes} : \C^T \sstimes \C^T \to \C^T$ such that the following square commutes:
$$
\begin{tikzcd}
  \C^T \sstimes \C^T
  \ar[r,"{\widetilde{\otimes}}"]
  \ar[d,"{U^T \sstimes U^T}"']
  &
  \C^T \ar[d,"{U^T}"]
  \\
  \C \sstimes \C
  \ar[r,"\otimes"']
  &
  \C
\end{tikzcd}
$$
that is, a lifting of $\otimes$. It is therefore encoded by a lifting datum $T_2 : T \circ \otimes \to (T \sstimes T) \circ \otimes$, which is exactly a colax magma structure $T_2$ such that $\mu$ and $\eta$ are colax magmatic.
\end{proof}

\subsection{Magma closed categories.} Let $\C$ be a magma category. Given $X$, $Y$ two objects of $\C$, a \emph{left internal Hom from $X$ to $Y$} is an object $\Lhom{X}{Y}$ of $\C$ which represents the presheaf
\begin{equation}
\label{eq:1}
  \pLhom{X}{Y}
  \colon \C^\opp \to\Set\qquad
  Z \mapsto \Hom(Z \otimes X,Y).
\end{equation}
In other words, it is an object $\Lhom{X}{Y}$ endowed with a morphism
$$e^X_Y : \Lhom{X}{Y} \otimes X \to Y,$$ the \emph{evaluation}, that enjoys the
following universal property: for all $f\colon Z\otimes X\to Y$ there exists a
unique $g\colon Z\to\Lhom{X}{Y}$ such that $e^X_Y(g\otimes 1)=f$. As a matter of course, the
representing object $\Lhom{X}{Y}$ is unique up to unique isomorphism.

The magma category $\C$ is \emph{left magma closed}, or simply \emph{left closed}, if any pair of objects $X$, $Y$ admits a left internal Hom $\Lhom{X}{Y}$. If such is the case, we have a functor
\begin{equation}
  \C^\opp \sstimes \C \to \Set\qquad
  (X,Y) \mapsto \Lhom{X}{Y}\label{eq:2}
\end{equation}
and for every object $X$ of $\C$, the endofunctor $\Lhom{X}{?}$ is right adjoint
to $? \otimes X$, with counit $e^X_Y : \Lhom{X}{Y} \otimes X \to Y$ and a unit
that we denote by $h^X_Y : Y \to \Lhom{X}{Y\otimes X}$.

\begin{rk}
If $\C$ is monoidal left closed, we have composition morphisms
$\circ_{X,Y,Z} : \Lhom{Y}{Z} \otimes \Lhom{X}{Y} \to \Lhom{X}{Z}$; this
composition for internal Homs is associative, with unit $h^X_\un : \un \to
\Lhom{X}{X \otimes \un} \cong \Lhom{X}{X}$. However, this is not the case for a general left closed magma category since the definition of composition involves the associativity constraint.
\end{rk}

One defines similarly a \emph{right internal Hom from $X$ to $Y$}, denoted by $\rhom{X}{Y}$, with its universal evaluation morphism
$X \otimes \rhom{X}{Y} \to Y$, and a \emph{right closed magma category}.

Since in this paper we consider only left internal Homs, from now on the superscript `$l$' will be omitted in the notation, so that we write $\plhom{X}{Y}$ and $\lhom{X}{Y}$ instead of $\pLhom{X}{Y}$ and $\Lhom{X}{Y}$.

\section{Slack closed functors and slack Hopf adjunctions}\label{sec:slack-closed}

\subsection{Slack closed functors}

 A functor between magma categories is slack left closed if it preserves left internal Homs in a `slack' way, that is, left internal Hom functors are preserved, but not the corresponding evaluations. The preservation is realised through an additional piece of data called a `slack left closed structure'. We give first a general definition, which does not assume the existence of internal Homs, and then we re-interpret this definition in special cases.

Let  $U \colon \D \to \C$ be a functor between magma categories.
A \emph{slack left closed structure} on $U$ is an isomorphism of presheaves, natural in $A$, $B$ objects of $\D$
\begin{equation}
  \label{eq:4}
  \mathfrak{E}^A_B\colon
 \int^{C\in \D} \C(-,U(C))\sstimes \plhom{A}{B}(C) \iso \plhom{U(A)}{U(B)}
\end{equation}
which means in particular that the coend on the left hand side exists.

Assume $U$ admits a cocontinuous extension $U_! \colon\widehat{\D} \to \widehat{\C}$ (see Section~\ref{sec:conventions}). This is true, for example, if $\D$ is essentially small or if $U$ has a left adjoint (which will always be the case in the rest of this article). Then a slack left closed structure on $U$ is just a natural isomorphism:
$$ \mathfrak{E}^A_B\colon U_!\plhom{A}{B}\to \plhom{U(A)}{U(B)}.$$

\emph{Slack right closed structures} are similarly defined.

\subsection{Slack preservation of internal Homs}\label{df:slack-preservation}

Let us now assume that the magma category $\D$ is left closed. One would expect
that a slack left closed functor $U\colon\D\to \C$ should preserve internal
Homs. The following definition makes this idea more precise. 

Let $U\colon\D\to \C$ be a functor from a left closed magma category $\D$ to a magma
category $\C$ and $\varphi_{A,B}\colon U(A)\otimes U(B)\to U(A\otimes B)$ be a
natural transformation (for $A, B$ in $\D$). We say that $(U,\varphi)$ \emph{slackly preserves left
internal Homs} if
  \begin{equation}
    \label{eq:evaluationphi}
    U(\lhom{A}{B})\otimes U(A)
    \xrightarrow{\varphi_{\lhom{A}{B},A}}
    U(\lhom{A}{B}\otimes A)\xrightarrow{U(e)}
    U(B)
  \end{equation}
is an evaluation\footnote{In the sense of left internal Homs.} in $\C$ for each evaluation $e\colon \lhom{A}{B}\otimes A\to B$ in $\D$.

The main point of introducing this new notion is that it encodes the slack preservation of left internal Homs independently of particular choices of internal Homs in $\D$. More precisely:

\begin{lemma}\label{lem-slack-preserve}
  Let $U \colon \D \to \C$ be a functor between magma categories, and assume $\D$ is
  left closed. Then there is a bijective correspondence between:
\begin{enumerate}
\item \label{item:1}
  slack left closed structures $\mathfrak{E}$ on $U$;
\item \label{item:2} given a choice $\lhom{\ }{\ }$ of left internal Homs in $\D$, families of transformations
  $$E = (E^B_A \colon{} U(\lhom{B}{A}) \otimes U(B) \to U(A))_{A,B \in \D}$$ natural in $A \in \D$ and dinatural in $B \in \D$ such that
  $(U(\lhom{B}{A}), E^B_A)$ is a left internal Hom from $U(B)$ to $U(A)$ in the
  category $\C$;
\item \label{item:3} natural transformations $\varphi$ such that $(U,\varphi)$ slackly preserves left internal Homs.
\end{enumerate}
The relation between the data $E$ and $\varphi$ is as follows:
$$\varphi_{A,B} = E^{B}_{A \otimes B} (U(h^{B}_A) \otimes UB), \quad E^B_A = U(e^B_A)\varphi_{\lhom{B}{A},B}\,.$$
\end{lemma}
\begin{proof}
Given a choice of left internal Homs in $\D$, the presheaf $\plhom{B}{A}$ is represented by  $\lhom{B}{A}$, so a slack left closed structure is an isomorphism of presheaves $\C(-,U\lhom{B}{A}) \cong \underline{\lhom{UB}{UA}}$, natural in $A, B$. By Yoneda
Lemma, this is precisely data $E$ as in part \ref{item:2} of the statement. For the bijection between \ref{item:2} and \ref{item:3}, by definition,  data $(U,\varphi)$ gives rise to data $E$ as prescribed in the statement. Conversely data $E$ gives rise to a natural transformation $\varphi$ as prescribed. The slack preservation condition on $\varphi$ is satisfied for $e=e^{B}_{A}$ since
\begin{align*}
U(e^B_{A})\varphi_{\lhom{B}{A},B}
&=
U(e^B_{A}) E^B_{\lhom{B}{A}\otimes B} (U(h^{B}_{\lhom{B}{A}}) \otimes UB) \\
&=E^B_{A}(U(\lhom{B}{e^B_{A}}  h^{B}_{\lhom{B}{A}}) \otimes UB) \blabla{(by naturality of $E$)} \\&= E^B_{A}.
\end{align*}
 That it is satisfied for all choices of $e$ results from the naturality of $\varphi$ and the uniqueness of left internal Homs up to unique isomorphism.
\end{proof}

\subsection{Slack Hopf adjunctions.}\label{sec-slack-Hopf-adj}

Let $\C$, $\D$ be magma categories, $F \adj U : \D \to \C$ be an adjunction with unit $\eta\colon \id_\C \rightarrow UF$ and counit
$\eps\colon FU \rightarrow \id_\D$. Define a \emph{slack left Frobenius
isomorphism} on $F \adj U$ to be a natural isomorphism
\begin{equation}
  \Psi_{X,A} : F(X \otimes U(A)) \to F(X) \otimes A \blabla{for $X \in \C$, $A \in \D$.}\label{eq:7}
\end{equation}
A \emph{slack left Hopf adjunction} is an adjunction between magma
categories equipped with a slack left Frobenius isomorphism.

\begin{thm}
  \label{thm-ll-fun-adj}
  Let $F\adj U$ be an adjunction between magma categories. There is a canonical
  bijection between:
\begin{enumerate}
\item slack left closed structures $\mathfrak{E}$ on $U$;
\item slack left Frobenius isomorphisms $\Psi$ on $F \adj U$.
\end{enumerate}
\end{thm}
\begin{proof}
  Since the functor $U$ has a left adjoint $F$, it has a cocontinuous extension $U_!: \widehat{\D} \to \widehat{\C}$ which coincides with $F^*$.
  A slack left closed structure on $U$ amounts therefore to an isomorphism
  $$F^*\plhom{A}{B}\iso\plhom{U(A)}{U(B)},$$ natural in $A,B$. The
  presheaf on the left hand side of the isomorphism is $\D(F(-)\otimes A,B)$,
  while that on the right is $\C(-\otimes U(A),U(B))$. The latter is isomorphic
  to $\D(F(-\otimes U(A)),B)$ by adjunction. It follows from Yoneda Lemma that  slack left
  closed structures on $U$ are in bijection with natural isomorphisms
  $F(-)\otimes A\cong F(-\otimes U(A))$.
\end{proof}

\begin{rk}\label{rk-phi-psi-mates} Let $F \adj U : \D \to \C$ be an adjunction between magma categories with $\D$ left closed. Then a slack left closed structure on $U$ corresponds bijectively with
\begin{enumerate}
\item a slack Frobenius isomorphism $\Psi$ (by Theorem~\ref{thm-ll-fun-adj});
\item a natural transformation $\varphi$ such that $(U,\varphi)$ slackly preserves left internal Homs (by Lemma~\ref{lem-slack-preserve}).
\end{enumerate}
It follows from the construction of these bijections that the natural transformation $\varphi$ is the mate of $\Psi$, that is:
\begin{equation*}
  \varphi_{A,B} =
  (U(\varepsilon_A\otimes B))(U\Psi_{U(A),B})
  \eta_{U(A)\otimes U(B)},
  \
  \Psi_{X,A}=
  \varepsilon_{FX\otimes A}
  F(\varphi_{FX,A})
  (\eta_X\otimes U(A)).
\end{equation*}
\end{rk}

\begin{rk}\label{rk-e-phi} If $\C$ and $\D$ are left magma closed, a slack left closed structure $\mathfrak{E}$ amounts to a natural isomorphism $\mathfrak{E}^B_A: U\lhom{B}{A} \iso \lhom{UB}{UA}$, and $E^B_{A} = e^{UB}_{UA} (\mathfrak{E}^B_A \otimes UB)$, hence explicit forms for the mates $\varphi$ and $\Psi$ are:
$$\varphi_{A,B} = e^{UB}_{U(A\otimes B)} (\mathfrak{E}^B_{A \otimes B} U(h^{B}_A)\otimes UB),$$
$$\Psi_{X,A} = \eps_{FX \otimes A} F\bigl( e^{UA}_{U(FX \otimes A)} (\mathfrak{E}^A_{FX\otimes A} U(h^A_{FX}) \eta_X \otimes UA)   \bigr) .$$
\end{rk}

\subsection{Slack Hopf colax magmatic adjunctions and slack Hopf structures.}

\begin{thm}\label{thm-ll-psi-beta}
Let $F \adj U: \D \to \C$ be a colax magmatic adjunction, and denote by $T=UF$ its colax magma monad. Then there is a canonical bijection:
$$\Theta^l:\Nat(F\otimes_\C(\id_\C \sstimes U), \otimes_{\D}(F\sstimes \id_\D)) \iso \Nat(\otimes_\C, \otimes_\C(T\sstimes T)),$$
defined by $\Psi \mapsto \Theta^l(\Psi) =\beta$, where:
$$
\left\{
\begin{array}{llll}
 & \beta_{X,Y} & = &U_2(FX,FY)\,U(\Psi_{X,FY}) \,\eta_{X \otimes TY}\,(X\otimes \eta_{Y}),\\
& \Psi_{X,A} & = & (FX \otimes \eps_A)\, \eps_{FX \otimes FUA}\, F(U^{-1}_2(FX,FUA))\, F(\beta_{X,UA}).\\
\end{array}
\right.$$

\end{thm}
\begin{proof} The bijection $\Theta^l$ results from the following sequence of natural bijections:
$\Nat(F\otimes_\C(\id_\C \sstimes U), \otimes_{\D}(F\sstimes \id_{\D})) \simeq
\Nat(\otimes_\C(\id_\C \sstimes U), U\otimes_{\D}(F\sstimes \id_{\D}))$ (by adjunction)
$ \simeq \Nat(\otimes_\C(\id_\C \sstimes U), \otimes_\C(T\sstimes U))$ (because $U$ is strong colax magmatic)
 $\simeq \Nat(\otimes_\C, \otimes_\C(T\sstimes T))$ again by adjunction. The formulae provided ensue.
\end{proof}

A \emph{slack left Hopf structure on the colax magmatic adjunction $F \adj U$}  with corresponding monad $T = UF$ is a natural transformation $\beta \in \Nat(\otimes_\C, \otimes_\C(T \sstimes T))$ such that $\Psi^\beta = {\Theta^l}^{-1}(\beta)$ is an isomorphism.

It follows tautologically  that we have a canonical bijection between slack left Hopf structures $\beta$ and left Hopf Frobenius isomorphisms $\Psi^\beta$.

\begin{rk}\label{rk:comp-lsh-adj} A composition of slack left closed functors being clearly slack left Hopf, the composition of slack left Hopf colax magmatic adjunctions
$F \adj U: \D \to \C$ and $G \adj V: \E \to \D$ is a slack left Hopf colax magmatic adjunction $G F\adj UV: \E \to \C$.

Slack left Hopf structures $\beta$, $\beta'$ for  $F \adj U$ and $G \adj V$ respectively give rise to a slack left Hopf structure $\beta''$ for
$GF \adj UV$ defined by
$$\beta''_{X,Y} = U_2(VGFX,VGFY)U(\beta'_{FX,FY})U^{-1}_2(FX,FY)\beta_{X,Y}.$$
\end{rk}

\section{Slack Hopf monads}
\label{sec:slack-hopf-adjunct}

\subsection{Slack Hopf monads and slack closed adjunctions.}

Let $T$ be a colax magma monad on a magma category $\C$. The \emph{left fusion operator of $T$} is the natural transformation
\begin{equation}
  H^l_{X,Y} \colon
  T(X\otimes TY)\xrightarrow{T_2(X,TY)}TX\otimes T^2Y \xrightarrow{1
    \otimes \mu_Y} TX\otimes TY
  \label{eq:8}
\end{equation}

For $\beta \in \Nat(\otimes_\C, \otimes_\C(T \sstimes T))$, define a natural transformation
$$\Phi^\beta_{X,Y} = (\mu_X \otimes \mu_Y) H^l_{TX,TY} T(\beta_{X,TY}) : T(X \otimes TY) \to TX \otimes TY.$$

A \emph{slack left Hopf structure on $T$} is a natural transformation $\beta$ such that $\Phi^\beta$ is an isomorphism.
Equivalently, it is a slack left Hopf structure on the adjunction $F^T\adj U^T$, but this is not obvious: see Theorem~\ref{thm-ll-psi-beta-2}.

A \emph{slack left Hopf  monad} is a colax magma monad equipped
with a slack left Hopf structure.

Slack left Hopf structures for a given colax magma monad are not unique, as a rule, but they form a torsor under certain group, see Corollary~\ref{th-uni}.

The morphism $\Phi^\beta$ can also be written as follows.
\begin{align*}
  \Phi^\beta_{X,Y} &
                     = (\mu_X \otimes \mu^2_Y)T_2(TX,T^2Y)T(\beta_{X,TY})\\
                   & = (\mu_X\otimes\mu_Y)T_2(TX,TY)T((TX\otimes\mu_Y)\beta_{X,TY}).
\end{align*}
The second expression means that $\Phi^\beta_{X,Y}$ is the morphism of
$T$-modules induced by $(TX\otimes\mu_Y)\beta_{X,TY}\colon X\otimes TY\to
TX\otimes TY$ under the free $T$-module adjunction.

\begin{exa}
  \label{ex:1} As one would expect, slack left Hopf monads generalise left Hopf monads in the sense of~\cite{MR2793022}. Indeed, a
  comonoidal monad $T$ on a monoidal category is left Hopf if and only if its left fusion
  operator~\eqref{eq:8} is invertible, which means exactly  that
  $\beta_{X,Y}=\eta_X\otimes\eta_Y$ is a slack left Hopf structure on $T$.
\end{exa}

\begin{rk}For a comonoidal monad $T$  on a monoidal category $\C$ to be a left Hopf monad depends therefore only on the magma structure of $\C$ and the colax magma monad structure of $T$, since the left fusion operator is defined only in these terms.\end{rk}

\begin{exa}\label{ex:alg-slack} Less expectedly, any monoid in a monoidal category $\C$ defines a slack left Hopf monad on $\C$. More precisely, let $(A,m,u)$ be a monoid in $\C$. Then the endofunctor $A \otimes ?$ of $\C$ has a structure of slack left Hopf monad on $\C$: the monad structure is the standard one, given by the product $m \otimes ?$ and unit $u \otimes ?$, the colax magma structure is $T_2(X,Y) = A \otimes X \otimes u \otimes Y$,
and $\beta = \eta \otimes \eta$ is a slack left Hopf structure, with corresponding $\Phi^\beta_{X,Y} = H^l_{X,Y} =\id_{A \otimes X \otimes A \otimes Y}$.
\end{exa}

\begin{thm}\label{thm-ll-psi-beta-2}
Let $F \adj U$ be a colax magmatic adjunction. Denote by $T=UF$ its colax magma monad. Then, for $\beta \in \Nat(\otimes_\C, \otimes_\C(T \sstimes T))$:
\begin{enumerate}
\item  $\Psi^\beta = {\Theta^l}^{-1}(\beta)$ and $\Phi^\beta$ are related by:
$$\Phi^\beta_{X,Y} =U_2(FX,FY) U(\Psi^\beta_{X,FY});$$
\item a slack left Hopf structure on $F\adj U$ is a slack left Hopf structure on $T$;
\item if $\Phi^\beta$ is an isomorphism, so is $U(\Psi^\beta)$, and we have
$$ U\left (({\Psi^\beta}_{X,A})^{-1}\right ) = T(X \otimes U\eps_{A})(\Phi^\beta_{X,UA})^{-1}(TX \otimes \eta_{UA})U_2(FX,A);$$
\item if $U$ is conservative,  a slack left Hopf structure on $T$ is a slack left Hopf structure on $F\adj U$.
\end{enumerate}
In particular the slack left Hopf structures of a colax magma monad $T$ coincide with those of the adjunction $F^T \adj U^T$.
\end{thm}

\begin{proof} Assertion (1) is straightforward. Assertion (2) means that given $\beta$, if  $\Psi^\beta$ is an isomorphism so is $\Phi^\beta$,
which follows from the fact that $U_2$ is an isomorphism. Assertion (4) is an immediate consequence of Assertion (3), which will result from the following
lemma.\phantom\qedhere
\end{proof}

\begin{lemma}
  \label{l:6}
  Let $G, H : \C \to \E$ be two functors, and consider a functor $U :
  \D\to \C$ with left adjoint $F$. If $f : G U \to H U$ is a natural
  transformation such that $f_{FU}$ is an isomorphism, then $f$ is an isomorphism and
  $f^{-1} = GU(\eps)f^{-1}_{FU} H(\eta_U)$.
\end{lemma}

\begin{proof} We will use the following classical fact.
Let $L\adj R : \mathcal{A} \to \mathcal{B}$ be an adjunction, with unit $\sigma\colon \id_\mathcal{A}\to RL$ and counit $\tau\colon LR\to \id_\mathcal{B}$.
Let $\mathcal{B}_0 \subset \mathcal{B}$ be the full subcategory of $\mathcal{B}$ defined by $\Ob(\mathcal{B}_0) = L(\Ob(\mathcal{A}))$.
Then $R_0 = R{\mid \mathcal{B}_0}$ is the right adjoint of the corestriction of $L$ to $\mathcal{B}_0$, and the counit of this adjunction is a split epimorphism
since $\tau_{LA}L\sigma_{A} = \id_{LA}$. As a consequence $R_0$ is
conservative (see the last two sentences in \cite[Section~3.4]{zbMATH02172008}).
If $f  \in \mathcal{B}(LA, LA')$ is such that $Rf$ is invertible, then $f$ is invertible and $f^{-1} = \tau_{LA} L(Rf)^{-1} L \sigma_{A'}$ (as a straightforward computation shows that this is a right inverse).

We apply the previous considerations to $\mathcal{A} = [\C,\E]$, $\mathcal{B} = [\D,\C]$. The functor $L = [U,\E] : \mathcal{A} \to \mathcal{B}$ has a right adjoint
$R = [F,\E]$. We obtain that if $f : L(G) = GU \to L(H) = HU$ is such that $R(F) = f_F$ is an isomorphism, then so is $f$ itself, and $f^{-1} = GU(\eps)f^{-1}_{FU} H(\eta_U)$.
\end{proof}

\begin{proof}[Proof of Theorem~\ref{thm-ll-psi-beta-2} cont.]
Now apply the lemma to $G, H : \C \to \C$, $G(Y) =T(X \otimes Y)$ and $H(Y)= TX \otimes Y$, and  $f : GU \to HU$,  $f_A = U_2(F(X), A)U(\Psi^\beta_{X,A})$. We have $f_{FY} = \Phi^\beta_{X,Y}$ so by the lemma, if $\Phi^\beta$ is invertible, so is $f$ for every choice of $X$, and $U(\Psi^\beta)$ as well since $U_2$ is an isomorphism. We have $ f^{-1}_A =  U(\Psi^\beta)^{-1}_{X,A} U_2(FX,A)^{-1}= T(X \otimes \eps_{UA}) \Psi^{-1}_{X,UA}(TX \otimes \eta_{UA})$, as expected, which concludes the proof of Assertion (3) and of the theorem.
\end{proof}

\begin{thm}\label{thm-ll-mon-mod} Let $T$ be a colax magma monad on a  magma category $\C$.  Then there are bijective correspondences between:
\begin{enumerate}[(i)]
\item  slack left Hopf structures $\beta$ on $T$;
\item slack left closed structures $\mathfrak{E}$ on $U^T$;
\item slack left Frobenius isomorphisms $\Psi$ on the adjunction $F^T \adj U^T$.
\end{enumerate}
The correspondence between $\beta$ and $\Psi$ is given by:
$$\beta_{X,Y} = (\Psi_{X,F^TY}) \eta_{X \otimes TY}(X\otimes \eta_Y) \blabla{and} \Psi_{X,(Y,\rho)}=(\mu_X \otimes \rho) H^l_{TX,Y} T(\beta_{X,Y}).$$
Moreover if $T$ is slack left Hopf and $\C$ is left closed, so is $\C^T$.
\end{thm}

\begin{proof} The bijection $\beta \leftrightarrow \Psi$ results directly from Theorem~\ref{thm-ll-psi-beta} and the fact that the forgetful functor of a monad is conservative. The formula for $\Psi$ results from that given in the theorem by straightforward computation.
The bijection $\Psi \leftrightarrow \mathfrak{E}$ is a special case of Theorem~\ref{thm-ll-fun-adj} and is only included for completeness' sake.
So the only really new statement is the last assertion.

Assume $\C$ is left closed and $T$ is slack left Hopf.
In order to show that $\C^T$ is left closed we use
\cite[Theorem~7.2]{zbMATH02222246} from Barr and Wells' book, which states that, given adjunctions $G \dashv V$ and $G' \dashv V'$, a functor $W$ such that $G'\cong WG$  has a right adjoint provided  $V$ is monadic and $W$ preserves $U$-split equalisers.
We apply Barr and Well's theorem to the adjunctions $F^T\dashv U^T$ and $F^T(-\otimes U^TA)\dashv [U^TA,-]U^T$ and the functor $-\otimes A$, where $A$ is an object of $\C^T$.
We dispose of the Frobenius isomorphism $F^T(- \otimes U^TA) \cong (- \otimes A)\, F^T$,  and the functor $- \otimes A$ preserves $U^T$-split coequalisers because
$U^T(-\otimes A)=U^T(-)\otimes U^T(A)$ preserves the said coequalisers and $U^T$ creates them. We conclude that $-\otimes A$ has a right adjoint for any $A$ in $\C^T$, that is, $\C^T$ is left closed.


\end{proof}

\subsection{Uniqueness of slack left Hopf structures.}\label{sec:unicity} Let $F\adj U$ be a colax magmatic adjunction, with $U$ conservative, and let $T = UF$ be its colax magma monad.

Observe that since, by definition, a slack left Frobenius isomorphism on $F\adj U$ is an isomorphism of functors from $F\otimes_\C(\id_{\C} \sstimes U)$ to $\otimes_{\D}(F\sstimes \id_\D)$, the group $\Aut(F\otimes_\C(\id \C \sstimes U))$ acts, on the right, transitively and freely on the set of slack left Frobenius isomorphisms. Now by Theorem~\ref{thm-ll-psi-beta} we have a canonical bijection
$$\Theta^l:\Nat(F\otimes_\C(\id_\C \sstimes U), \otimes_{\D}(F\sstimes \id_\D)) \iso \Nat(\otimes_\C, \otimes_\C(T\sstimes T))$$
in which slack left Frobenius isomorphisms correspond bijectively with slack left Hopf structures on $T$.

On the other hand, we have also by adjunction a canonical bijection
$$\Xi^l:\End(F\otimes_\C(\id_\C \sstimes U)) \iso \Nat(\otimes_\C, T\otimes_\C(\id_\C \times T))$$
which allows us to endow $\Nat(\otimes_\C, T\otimes_\C(\id_\C \times T))$ with a monoid product, denoted by $*$, in such a way that $\Xi^l$
becomes an isomorphism of monoids.

Denote by $\M^l(T)$ the monoid  $\Nat(\otimes_\C, T\otimes_\C(\id_\C \times T))$ and by $\G^l(T)$ its group of invertibles. It then follows immediately from the above considerations that $\G^l(T)$ acts (on the right) transitively and freely on the set of slack left Hopf structures of $T$. Besides, this is true for any colax magma monad $T$, since the adjunction $F^T \adj U^T$ is colax magmatic and $U^T$ is conservative.

A natural endomorphism $\Gamma_{X,A}$ of $F(X \otimes UA)$ and the corresponding natural transformation
$\gamma_{X,Y} : X \otimes Y \to T(X \otimes TY)$ are related by the following formulae:
$$\gamma_{X,Y} = U(\Gamma_{X,FY}) \eta_{X \otimes TY}(X \otimes \eta_Y) \ \text{and}\ \Gamma_{X,A} = F(X \otimes U\eps_A)\eps_{F(X \otimes TUA)} F(\gamma_{X,UA}),$$
which allows us to compute the convolution product of $\M^l(T)$, given by: 
$$(\gamma*\gamma')_{X,Y} = \mu_{X \otimes TY}T^2(X \otimes \mu_Y)T(\gamma_{X,TY})\gamma'_{X,Y},$$
and the action of  $\M^l(T)$  on the right on $\Nat(\otimes_\C, \otimes_\C(T\sstimes T))$, given  by:
$$(\beta \lhd \gamma)_{X,Y} = (\mu_X \otimes \mu_YT(\mu_Y)) T_2(TX,T^2Y)T(\beta_{X,TY})\gamma_{X,Y}
= \Phi^\beta_{X,Y} \gamma_{X,Y}.$$
This is summed up in the following theorem.

\begin{thm}\label{th-uni} Let $T$ be a slack left Hopf monad. Then  the set of slack left Hopf structures on $T$ is a torsor under the right action of the group $\G^l(T)$.
\end{thm}

\begin{cor} \label{cor-slh-h}
Let $T$ be a slack left Hopf monad on a magma category $\C$, and let $\beta$ be a slack left Hopf structure on $T$. Under the above notations, define  $\alpha \in \M^l(T)$ by
$$\alpha_{X,Y} = ({\Phi^\beta_{X,Y}})^{-1} (\eta_X \otimes \eta_Y).$$
Then the following conditions are equivalent:
\begin{enumerate}[(i)]
\item The left fusion operator $H^l$ is an isomorphism;
\item $\eta \otimes \eta$ is a slack left Hopf structure on $T$;
\item $\alpha \in \G^l(T)$.
\end{enumerate}
\end{cor}

\begin{proof}
Clearly (i) and (ii) are equivalent, since $H^l = \Phi^{\eta \otimes \eta}$.  Now by  Theorem~\ref{th-uni}, since $\beta$ is a slack left Hopf structure, (ii) is equivalent to the existence of $\gamma \in \G^l(T)$
such that $\beta \lhd \gamma = \Phi^\beta \gamma = \eta \otimes \eta$, but such a $\gamma$ must be equal to $\alpha$, since $\Phi^\beta \alpha = \eta \otimes \eta$ and $\Phi^\beta$ is an isomorphism.  This shows the equivalence of (ii) and (iii).
\end{proof}

\subsection{Slack antipodes} Let $T$ be a slack left Hopf monad on a left closed magma category $\C$, and let $\beta$ be a slack left Hopf structure on $T$. It results from Theorem~\ref{thm-ll-mon-mod} that $\C^T$ is left closed and the forgetful functor $U$ preserves left internal Homs (but not necessarily their adjunction units and counits).
In this section we make this explicit, in particular we give a notion of slack left antipode associated with a left Hopf structure and express the action of $T$ on the left internal Hom between $T$-modules in terms of the slack left antipode.

\begin{thm}\label{thm-antipode} Let $T$ be a slack left Hopf monad on a left closed magma category $\C$, with slack left Hopf structure $\beta$. For $X,Y$ in $\C$ set
$$\xi_{X,Y}= T(X \otimes \mu_Y)(\Phi^{\beta}_{X,TY})^{-1}(\mu_X \otimes T^2Y)\beta_{TX,TY} : TX \otimes TY \to T(X \otimes TY)$$ and define the \emph{slack left antipode associated to $\beta$}, $S^l$, by
$$S^l_{X,Y}: = \lhom{Y}{Te^{TY}_X} \lhom{\eta_Y}{\xi_{\lhom{TY}{X},Y}} h^{TY}_{T\lhom{TY}{X}}:T\lhom{TY}{X} \to \lhom{Y}{TX}$$

Given two $T$-modules $(X,r)$ and $(Y,\rho)$, $\lhom{Y}{X}$ is a $T$-module with action $R: T(\lhom{Y}{X}) \to \lhom{Y}{X}$ given by:
$$R=\lhom{Y}{r} S^l_{X,Y} T\lhom{\rho}{X}.$$
The $T$-module so defined is a left internal Hom from $(Y,\rho)$ to $(X,r)$, with universal evaluation (i.e. adjunction counit)
$e^{(Y,\rho)}_{(X,r)} : \lhom{Y}{X} \otimes Y \to X$  defined by:
$$e^{(Y,\rho)}_{(X,r)}=r T(e^Y_X)T(\lhom{Y}{X} \otimes \rho)\alpha_{\lhom{Y}{X},Y},$$
where $\alpha = \Phi^{\beta -1}(\eta \otimes \eta)$, the adjunction unit being given by:
$$h^{(Y,\rho)}_{(X,r)}= \lhom{Y}{(r \otimes \rho)\beta_{X,Y}} h^Y_X.$$
\end{thm}
\begin{proof} Observe first that $\C^T$ is left closed by Theorem~\ref{thm-ll-mon-mod}, and  $\beta$ defines a left closed structure $\mathfrak{E}^{(Y,\rho)}_{(X,r)} : U(\lhom{(Y,\rho)}{(X,r)}) \iso \lhom{Y}{X}$ for $(X,r)$, $(Y,\rho)$, in $\C^T$. By transport of structure we may choose left internal Homs in $\C^T$ in such a way that $\mathfrak{E}$ is an identity.

Fix $\widetilde{Y} = (Y,\rho)$ once and for all, and consider the adjunction $G \adj V$, where $G = ? \otimes Y$ and $V= \lhom{Y}{?}$.
The functor $G : \C \to \C$ lifts to a functor $\widetilde{G} = ? \otimes \widetilde{Y} : \C^T \to \C^T$, with corresponding distributive law $\ZETA : TG \to GT$  given by:
$$\ZETA_X = (TX \otimes \rho)T_2(X,Y): T(X \otimes Y) \to X \otimes T(Y).$$
Let $\widetilde{V} = \lhom{\widetilde{Y}}{?} :\C^T \to \C^T$, so that we have an adjunction $\widetilde{G}\adj \widetilde{V}$. With our choices of left internal Homs in $\C^T$, $U\widetilde{V} = VU$, that is, $\widetilde V$ is a lift of $V$.

The following lemma will allow us to better understand the adjunction $\widetilde G \adj \widetilde V$  in terms of $\ZETA$ and $\BETA$.

\begin{lemma}
Let $G\adj V : \C\to \D$ be an adjunction, with unit $h: \id_\C \to VG$ and counit $e: GV \to \id_\D$. Let $(T,\mu,\eta)$ be a monad on $\C$, $(P,m,u)$ be a monad on $\D$, and let $\widetilde{G} \adj \widetilde{V} : \C^T \to \D^P$ be an adjunction, with unit $\widetilde{h}$ and counit $\widetilde{e}$, such that $\widetilde{G}$ is a lift of $G$ and $\widetilde{V}$ is a lift of $V$.

Denote by $\ZETA : PG \to GT$ the lifting datum for $\widetilde G$,
and by $\XI: GT \to PG$ the mate of the lifting datum $TV \to VT$ for $\widetilde V$. Let $\BETA: G \to  GT$ be the mate of the lifting datum $\id_\C \to VGT$ of $\widetilde{h}$ and $\ALPHA: G\to PG$ the mate of the lifting datum $GV \to P$ of $\widetilde e$.

Then the following identities hold:

$$\begin{array}{ll}
(a)\quad \ZETA m_G = G(\mu)\ZETA_T T(\ZETA)\quad &(a')\quad  \ZETA u_G = G(\eta)\\
(b) \quad \XI G(\mu) = m_G P(\XI)\XI_T\quad &(b') \quad \XI_ G(\eta) = u_G\\
(c) \quad m_G P(\XI) \ALPHA_T\ZETA = m_G P(\ALPHA)\quad &(c')\quad m_G P(\XI)\ALPHA_T \BETA = u_G\quad \\
(d) \quad G(\mu) \ZETA_T P(\BETA) \XI = G(\mu)\BETA_T\quad &(d')\quad G(\mu)\ZETA_T P(\BETA)\ALPHA = G(\eta)
\end{array}$$

Moreover $\XI$, $\ALPHA$, $\BETA$ are determined by $\PSI = G(\mu)\ZETA_TP(\BETA): PG \to GT$. More precisely:
\begin{enumerate}
\item $\PSI$ is invertible, and  $\PSI^{-1} = m_G P(\XI) \ALPHA_T$;
\item $\BETA = \PSI u_G$;
\item $\ALPHA = \PSI^{-1}G(\eta)$;
\item $\XI = \PSI^{-1}G(\mu)\BETA_T$.
\end{enumerate}
\end{lemma}

\begin{proof} The identities (a) --- (d')  translate what we know about $\widetilde{G} \adj \widetilde{V}$: $\ZETA$ is a lifting datum for $G$ ((a), (a')), $\XI$ mate of a lifting datum for $V$ ((b), (b')),
$\ALPHA$ encodes $\widetilde e : \widetilde{G}\widetilde{V} \to \id_{\C^T}$ (c), $\BETA$ encodes $\widetilde h: \id_{\D^P} \to \widetilde{V}\widetilde{G}$ (d), and the adjunction axioms are (c') and (d').

Let $\THETA = m_GP(\ZETA)\ALPHA_T$.
Then $\PSI\THETA = G(\mu)\ZETA_TP(\BETA)m_GP(\XI)\ALPHA_T= \id_{GT}$ (using (a), (d) and (d')), and similarly $\THETA\PSI = \id{PG}$ (using (b), (c) and (c')), hence Assertion (1).
Also $\PSI u_G =G(\mu)\ZETA_TP(\BETA)u_G = G(\mu)\ZETA_T u_{GT}\BETA = G(\mu)G(\eta_T)\BETA$ (using (b')) $ = \BETA$  and similarly $\THETA G(\eta) = \ALPHA$ (using (a')), hence Assertions (2) and (3).
Lastly $\THETA G(\mu)\BETA_T = m_G P(\XI)\ALPHA_T G(\mu)\BETA_T = \XI$ (using (b) and (c)), hence Assertion (4).
\end{proof}

\begin{rk} The conditions (a) --- (d') are not only necessary but also sufficient for the data $(\ZETA, \XI,\ALPHA, \BETA)$ to define a lifted adjunction $\widetilde{G} \adj \widetilde{V}$ satisfying the hypotheses of the lemma, as will be shown in a forthcoming paper. 
\end{rk}

Let us return to Theorem \ref{thm-antipode} and apply this lemma with $\D =\C$ and $P=T$. We first show  $$\BETA_X = (TX \otimes \rho) \beta_{X,Y} \blabla{and}
\PSI_X = U\Phi^\beta_{X,\widetilde{Y}}.$$
Since $\BETA: G \to GT$ is the mate of the lifting datum of $\widetilde{h}$, we compute:
$$\BETA_X =e^Y_{TX \otimes Y}(U(h^{\widetilde{Y}}_{FX} ) \eta_X \otimes Y).$$
In view of Remark~\ref{rk-e-phi}, $\mathfrak{E}$ being an identity, the left Frobenius isomorphism $\Psi^\beta$ is: $$\Psi_{X,\widetilde{Y}} = \eps_{FX \otimes \widetilde{Y}} F\bigl( e^Y_{TX \otimes Y}(U(h^{\widetilde{Y}}_{FX})\eta_X \otimes Y \bigr) = \eps_{FX \otimes \widetilde{Y}}F(\BETA_X).$$
So $U\Psi^\beta_{X,\widetilde{Y}} \eta_{X \otimes Y} = U\eps_{FX \otimes \widetilde{Y}} T\bigl(\BETA_X) \eta_{X \otimes Y} = U\eps_{FX \otimes \widetilde{Y}} \eta_{TX \otimes Y} \BETA_X = \BETA_X$. On the other hand
$U\Psi^\beta_{X,\widetilde{Y}}\eta_{X \otimes Y} =(\mu_X \otimes \rho\mu_Y)T_2(TX,TY) T(\beta_{X,Y})\eta_{X \otimes Y} =  (TX \otimes \rho)\beta_{X,Y}$, hence  $\BETA_X = (TX \otimes \rho)\beta_{X,Y}$. Now
$$\PSI_X = (\mu_X \otimes \rho)T_2(TX,Y)   T((TX \otimes \rho)\beta_{X,Y}) = U\Psi^\beta_{X,\widetilde{Y}}.$$
We can now express $\PSI$ and $\PSI^{-1}$ in terms of $\Phi$ (the latter, using Theorem~\ref{thm-ll-psi-beta-2}):
 $$\PSI_X = (TX \otimes \rho)\Phi_{X,Y}T(X \otimes \eta_Y) \blabla{and}\PSI^{-1}_X = T(X \otimes \rho) \Phi^{-1}_{X,Y}(TX \otimes \eta_Y).$$
Having identified $\ZETA$, $\BETA$, $\PSI$ and $\PSI^{-1}$, we obtain $\ALPHA_X = \PSI^{-1}_X(\eta_X \otimes Y)$ and $\XI_X = \PSI^{-1}_X(\mu_X \otimes Y)\BETA_{TX}$.

We have $\ALPHA_X = (TX \otimes \rho) \Phi^{-1}_{X,Y}(\eta_X\otimes \eta_Y) : G(X) \to TG(X)$, its mate $\alpha^* : GV(X) \to T(X)$ is:
$\ALPHA^*_X = T(e^Y_X)\ALPHA_{\lhom{Y}{X}}$ and so $$\widetilde{e}^{\widetilde{Y}}_{(X,r)} = r\ALPHA^*_X =
rT(e^Y_X)(T\lhom{Y}{X} \otimes \rho) \Phi^{-1}_{\lhom{Y}{X},Y}(\eta_{\lhom{Y}{X}}\otimes \eta_Y).$$
Similarly $\BETA_X = (TX \otimes \rho) \beta_{X,Y}$ gives us $h^{\widetilde{Y}}_{(X,r)}= \lhom{Y}{(r \otimes \rho)\beta_{X,Y}} h^Y_X$.

Lastly, we compute $\XI = \PSI^{-1}G(\mu)\BETA_T: GT \to TG$:
\begin{align*}\XI_X &=  T(X \otimes \rho) \Phi^{-1}_{X,Y}(\mu_X \otimes \eta_Y\rho)\beta_{TX,Y}
= T(X \otimes \rho) \xi_{X,Y}(\mu_X \otimes \eta_Y)\beta_{TX,Y}.
\end{align*}
The lifting datum $\XI^*: TV \to VT$ for $\tilde{V}$ is the mate of $\XI$, that is:
\begin{align*}
\XI^*_X &= (T(e^Y_X) \otimes Y)\lhom{X}{\XI_{\lhom{Y}{X}}}h^Y_{T\lhom{Y}{X}} \\&=
(T(e^Y_X) \otimes Y)\lhom{X}{
(T(\lhom{Y}{X} \otimes \rho) \xi_{\lhom{Y}{X},Y}(\lhom{Y}{X} \otimes \eta_Y)
}h^Y_{T\lhom{Y}{X}} = S^l_{X,Y} T\lhom{\rho}{X},
\end{align*}
so we get $R =\lhom{Y}{r}\XI^*_X = \lhom{Y}{r} S^l_{X,Y} T\lhom{\rho}{X}$.
\end{proof}

\begin{exa}Given a magma category $\mathcal{C}$, a slack left Hopf structure on the
  identity functor $1_{\mathcal{C}}$ (taken as a trivial colax magma monad) amounts to a natural isomorphism
  \begin{equation}
    X\otimes Y\cong X\otimes Y.
  \end{equation}
On the category $\mathbf{Ab}$ of abelian
  groups, we have a natural isomorphism as above given by $x\otimes y\mapsto
  -x\otimes y$. This corresponds to the left closed structure with the usual
  internal Homs but evaluation given by $[X,Y]\otimes X\to Y$, $f\otimes x\mapsto
  -f(x)$.
  \end{exa}

\subsection{Push-forwards and cross-products of slack Hopf monads}
    Consider an adjunction $F \adj U: \D\to \C$ and a monad $Q$ on $\D$. Then $UQF$ has a structure of a monad on $\C$, being the monad of the composition $F^QF \adj UU^Q:\D^Q\to \C$ of the adjunctions $F \adj U$ and $F^Q \adj U^Q$. This monad is called the \emph{push forward of $Q$ along $F \adj U$.}

If $F\adj U$ is a colax magmatic adjunction between magma categories and if $Q$ is a colax magma monad on $\D$, then its push forward along $F \adj U$ is a colax magma monad on $\C$.
Moreover, by virtue of Remark~\ref{rk:comp-lsh-adj}, if $Q$ is a slack left Hopf monad on $\D$, and $F \adj U: \D \to \C$ is a slack left Hopf adjunction, then the push forward of $Q$ along $F\adj U$ is a slack left Hopf monad on $\C$.

Let $T$ be a monad on a category $\C$. If $Q$ is a monad on the category $\C^T$ of $T$\ti modules,
the monad of the composite adjunction  $$F^Q F^T \adj U^TU^Q: (\C^T)^Q \to \C,$$ that is, the push forward of $Q$ along $F^T \adj U^T$, is called the \emph{cross product of $T$ by $Q$} and denoted by $Q \cp T$ (see \cite[Section 3.7]{BV3}).  As an endofunctor of~$\C$, $Q \cp T=U^TQF^T$. The product $p$ and unit $e$ of $Q \cp T$ are:
\begin{equation*}
p=q_{F^T} Q(\varepsilon_{Q F^T})  \quad \text{and}\quad e=v_{F^T}\eta,
\end{equation*}
where $q$ and $v$ are the product and the unit of $Q$, and $\eta$ and $\varepsilon$ are the unit and counit of the adjunction $(F^T,U^T)$.

If $\C$ is a magma category, $T$ is a colax magma monad on $\C$ and $P$ a colax magma monad on $\C^T$, then
$Q \cp T$ is a colax magma monad on $\C$.

If in addition $T$ and $Q$ are slack left Hopf monads, so is $Q \cp T$, and slack left Hopf structures  $\beta$ and $\beta'$ for $T$ and $Q$ respectively induce a slack left Hopf structure $\beta''$ on $Q\cp T$ given by:
$$\beta''_{X,Y} = \beta'_{F^TX,F^TY} \beta_{X,Y}.$$


%

\section{The cartesian case}\label{sec-unit-cart}
\subsection{Slack Hopf monads on cartesian categories}
\label{sec:slack-cartesian-categories} If $\C$ is a cartesian category, that is, a category having finite products, we denote by $\pt$ its terminal object by $\pt_X$ the unique morphism $X \to \pt$ for $X$ object of $\C$.

A cartesian category $\C$ will be viewed as a monoidal category $(\C, \sstimes, \pt)$.
Any functor $F : \C \to \D$ between cartesian categories admits a unique comonoidal structure $(F,F_2,F_0)$, with
$$F_2(X,Y)  = (F(X \sstimes \pt_Y),F(\pt_X \sstimes Y)) : F(X \sstimes Y) \to FX \sstimes FY \blabla{and} F_0 = \pt_{F(\pt)}.$$
Moreover, any natural transformation between functors between cartesian categories is comonoidal.
Consequently a monad on a cartesian category becomes a comonoidal monad when equipped with its unique comonoidal structure.

\begin{thm}\label{th-hopf-cart} Let $T$ be a monad on a cartesian category, endowed with its unique comonoidal structure. Then $T$ is slack left Hopf if and only if it is left Hopf.
\end{thm}

We will prove this Theorem in Section~\ref{sec:proof-cart}, but first we need to introduce the notions of unitary magma categories and unitary colax magma monads.

\subsection{Slack Hopf monads in the presence of a unit object}
Let $\C$ be a magma category. A \emph{left} (resp. \emph{right}) \emph{unit object} for $\C$ is an object $\un$ of $\C$ endowed with a natural isomorphism
$$l_X \colon \un \otimes X \to X \blabla{(resp. $r_X \colon X \otimes \un \to X$)}\blabla{for} X \in \C.$$ A \emph{left} (resp. \emph{right}) \emph{unitary magma category} is a magma category equipped with a left (resp. right) unit object.

Given a left (resp. right) unitary magma category $\C$ with left unit object $(\un,l)$ (resp. right unit object $(\un,r)$), a \emph{right} (resp. \emph{left}) \emph{unitary colax magma monad on $\C$} is a colax magma monad $T$ on $\C$ endowed with a morphism $T_0 : T\un \to \un$ satisfying:
\begin{enumerate}
\item  $T_0\mu_\un = T_0 T(T_0)$ and $T_0 \eta_{\un} = \id_\un$;
\item for $X$ in $\C$, $(T_0 \otimes TX)T_2(\un,X) = l^{-1}_{TX}T(l_X)$ (resp. $(TX \otimes T_0)T_2(X,\un) = r^{-1}_{TX}T(r_X)$).
\end{enumerate}
The first condition means that $\widetilde{\un} = (\un,T_0)$ is a $T$-module, and the second, that $l_{U^T}$ (resp. $r_{U^T}$) is $T$-linear, making $\widetilde{\un}$ into a left (resp right) unit object in $\mathcal{C}^T$.

\begin{rk} A monoidal category $\C$ is left and right unitary as a  magma category, and a comonoidal monad on it is left and right unitary as a colax magma monad.
\end{rk}
We will need two rather technical lemmas, the first in the right unitary context, the second in the left unitary context, both of which will be used in the cartesian case.

Let $T$ be a right unitary colax magma monad on a right unitary magma category $\C$. Given a natural transformation $\beta\colon \id_\C \otimes \id_\C \to T \otimes T$, we define a natural transformation
$\beta^u \colon \id_\C \to T$ by setting:
\begin{equation*}
    \beta^u_X\colon
    X\xrightarrow{r_X^{-1}}X\otimes \un\xrightarrow{\beta_{X,\un}}
    TX\otimes T\un
    \xrightarrow{TX\otimes T_0}
    TX\otimes \un
    \xrightarrow{r_{TX}}
    TX
\end{equation*}
We say that a slack left Hopf structure $\beta$ on $T$ is \emph{normalised} if $\beta^u = \eta$.

\begin{lemma} \label{lem-unitary-normalised}
Let $T$ be a right unitary colax magma monad on a right unitary magma category $\C$. Let $\beta$ be a slack left Hopf structure on $T$. Then
\begin{enumerate}
\item \label{item:betau1} $\beta^u$ is convolution invertible,

\item \label{item:betau2}
    denoting by $u$ the convolution inverse of $\beta^u$, the natural transformation $\overline{\beta}$ defined by $\overline{\beta}_{X,Y} = (\mu_X Tu_X \otimes TY)\beta_{X,Y}$ is a normalised slack left Hopf structure on $T$.
\end{enumerate}
\end{lemma}

\begin{proof}
For $X$ in $\C$, $B_X =
r_{TX}\Psi^\beta_{X,\widetilde{\un}} T(r^{-1}_X)$, being a natural automorphism of $F^T(X)$, is therefore of the form $\hat{b}_X$, where $b$ is an convolution invertible element of $\Nat(\id_\C,T)$ (see Section~\ref{sec:convo}). We have
\begin{multline*}
b_X= B_X \eta_X = r_{TX}(\mu_X \otimes T_0\mu_\un)T_2(TX,T\un)T(\beta_{X,\un})T(r^{-1}_X)\eta_X \\
= r_{TX}(\mu_X \otimes T_0\mu_\un)T_2(TX,T\un)\eta_{TX \otimes T\un}\beta_{X,\un}\,r^{-1}_X
= r_{TX}(TX \otimes T_0)\beta_{X,\un}r^{-1}_X = \beta^u_X,
\end{multline*}
so $\beta^u$ is convolution invertible, which proves Assertion (\ref{item:betau1}).

Now let $\overline{\Phi} = \Phi^{\overline{\beta}}$, with $\overline{\beta}_{X,Y} = (\mu_X Tu_X \otimes TY)\beta_{X,Y}$. A straightforward computation shows that
$\overline{\Phi}_{X,Y} = (\mu_X Tu_X \otimes TY) \Phi_{X,Y}$, and since $\mu_X \, Tu_X$ and $\Phi_{X,Y}$ are isomorphisms, so is $\overline{\Phi}_{X,Y}$; in other words, $\overline{\beta}$ is a slack left Hopf structure.
Moreover,
\begin{align*}\overline{\beta}^u_X & =r_{TX} (\mu_X Tu_X \otimes \un)\tau_X r^{-1}_X= \mu_X T(u_X) r_{TX} \tau_X r^{-1}_X
 = \mu_X T(u_X) \beta^u_X = \eta_X,
\end{align*}
that is, $\overline{\beta}$ is normalised, which proves Assertion (2).
\end{proof}

\begin{lemma} \label{lemma-unitary-delta}
Let $T$ be a left unitary colax magma monad on a left unitary magma category $\C$. Assume $T$ has a slack left Hopf structure on $T$ of the form $\beta = \eta \otimes \delta$,
where $\delta$ is a natural transformation  $\id_\C \to T$. Then
\begin{enumerate}
\item The natural transformation $\delta$ has a left inverse  for the convolution product. Such a left inverse $\theta$ is given by $\theta_Y = \mu_Y T(l_{TY}) \alpha_{\un,Y} l^{-1}_Y$.
\item If $\delta$ has a left inverse which is colax magmatic as a natural transformation $\id_\C \to T$, then the fusion operator $H^l$ is an isomorphism.
\end{enumerate}
\end{lemma}

\begin{proof} Let $\beta = \eta \otimes \delta$.
We have $\Phi^\beta_{X,Y} \eta_{X\otimes TY}(X \otimes \eta_Y) = \beta_{X,Y} = \eta_X \otimes \delta_Y$.
On the other hand $\Phi^\beta\alpha = \eta \otimes \eta$, from which results
$$\Phi^\beta_{X,Y} T(X \otimes \mu_Y T\delta_Y)\alpha_{X,Y} = (TX \otimes \mu_Y T\delta_Y)\Phi^\beta_{X,Y}\alpha_{X,Y} = \eta_X \otimes \delta_Y,$$
hence by monicity of $\Phi^\beta$: $T(X \otimes \mu_Y T\delta_Y) \alpha_{X,Y} = \eta_{X \otimes TY}(X \otimes \eta_Y)$.
Applying this to $X = \un$, precomposing with $l^{-1}_Y$ and postcomposing with $\mu_Y T(l_{TY})$, we obtain
$$\mu^2_Y T^2\delta_Y T(l_{TY})\alpha_{\un,Y}l^{-1}_Y = \mu_Y T\delta_{Y}\mu_Y T(l_{TY})\alpha_{\un,Y}l^{-1}_Y = (\theta*\delta)_Y) \eta_Y$$
where $\theta_Y = \mu_Y T(l_{TY})\alpha_{\un,Y}l^{-1}_Y$,  that is Assertion (1).

On the other hand, setting  $\check{\delta} = U(\eps) \delta_U$, we have
$$\Psi_{X,A} = \Psi^0_{X,A}F(X \otimes \check{\delta}_A),$$
where $\Psi^0$ is the fusion operator $(FX\otimes \eps_A)F_2(X,UA)$. In particular, since $\Psi$ is an iso $\Psi^0$ is split epic.
%
%
Now assume $\theta$ is a colax magmatic left inverse of $\delta$ (not necessarily of the form given above).
Let $\hat\theta = \eps_F F(\theta)$. In the monoid isomorphism $\Nat(\id_\C,T)^\opp \iso \Hom(F,F)$, $\theta$ is sent to $\hat\theta$, so $\hat\theta$ is a split monomorphism.
We claim that the following identity holds:
$$\Psi_{X,A} \hat\theta_{X \otimes UA} = (\hat\theta_X \otimes A) \Psi^0_{X,A}.$$

Before proceeding to prove this identity, let us show how it allows us to conclude. We already know that $\Psi^0$ is split epic, and this formula shows that it is monic because $\Psi$ is an isomorphism and $\hat\theta$ is split monic, therefore $\Psi^0$ is an isomorphism, and so is $H^l_{X,Y} = U(\Psi^0_{X, FY})$.

Since $\theta$ is colax magmatic, so is $\hat\theta$. We compute:
\begin{align*}
\Psi_{X,A}&\hat\theta_{X \otimes UA}  = \Psi^0 F(X \otimes \check\delta_A)\hat\theta_{X \otimes UA}
= (FX \otimes \eps_A)F_2(X,UA)\hat\theta_{X \otimes UA}F(X \otimes \check\delta_A)\\
&= (\hat\theta_X \otimes \eps_A\hat\theta_{UA})F_2(X,UA)F(X \otimes \check\delta_A)
= (\hat\theta_X \otimes \eps_A\hat\theta_{UA} F\check\delta_A)F_2(X,UA).
\end{align*}

Finally, we have
\begin{align*}\eps\hat\theta_U F \check\delta & = \eps(\eps FU)(F\theta U)(FU\eps)(F\delta U) =\eps(FU\eps)(F\theta U)(FU\eps)(F\delta U)\\
& =\eps(F \check\theta)(F\check\delta) = \eps,
\end{align*}
hence $\Psi_{X,A}\hat\theta_{X \otimes UA} = (\hat\theta_X\otimes A)\Psi^0_{X,A}$
and we are done.
\end{proof}

\subsection{Proof of Theorem~\ref{th-hopf-cart}}\label{sec:proof-cart} Let $T$ be a (comonoidal) monad on a cartesian category $\C$.
If $T$ is left Hopf, it is slack left Hopf with slack left Hopf structure $\beta = \eta \sstimes \eta$.
Conversely, assume $T$ is slack left Hopf, with slack left Hopf structure $\beta$.
Since $\C$ is cartesian, we have
$\beta = \gamma \sstimes \delta$, where $\gamma, \delta  \in \Nat(\id_\C, T)$.
Since $\C$ is monoidal and $T$ is comonoidal, Lemma~\ref{lem-unitary-normalised} applies, and we may assume that $\beta$ is normalised, that is $\beta^u = \gamma = \eta$, and so
$\beta = \eta \sstimes \delta$. We may therefore also apply Lemma~\ref{lemma-unitary-delta}, which tells us (1) that $\delta$ has a left inverse $\theta$ for the convolution product: $\mu T(\delta) \theta = \eta$,
 and (2) that $H^l$ is an isomorphism, since $\theta$ is automatically comonoidal.
\qed

\goodbreak
\subsection{Groupoids as slack Hopf categories}\label{sec:cat}

\subsubsection{Categories as colax magma monads}
Let $\A$ be a small category, and denote by $S = |\A|$ its set of objects, viewed as a discrete category. The inclusion functor $J : S \to \A$ induces by restriction a forgetful functor
$U = [\A,\Set] \to [S,\Set]$ which is monadic. The induced monad $T$ on $[S,\Set]$ sends $X = (X_s)_{s\in S} \in [S,\Set]$ to $T(X) \in [S,\Set]$ defined by
$T(X)_s = \sscoprod_{t \in S} \A(t,s) \sstimes X_t$ for $s \in S$. The monad $T$ is cocontinuous because so is $U$.
The product $\mu$ and unit $\eta$ of $T$ are as follows:
\begin{align*}(\mu_X)_s &: \sscoprod_{t,t' \in S}  \A(t',s) \sstimes \A(t, t') \sstimes X_t \to \sscoprod_{t \in S} \A(t,s) \sstimes X_t, \quad (f,g,x) \mapsto (f\circ g,x), \\
(\eta_X)_s &: X_s \mapsto \sscoprod_{t \in S}  \A(t,s) \sstimes X_t, \quad x \mapsto (\id_X,x).
\end{align*}
We view $[S,\Set]$ as a monoidal category for the cartesian product. The unique comonoidal structure on $T$  is given by $T_0:  \left (\sscoprod_{t\in S}  \A(t,s)\right )_{s\in S}\to (\pt_s)_{s\in S}$ the unique possible morphism (that is $T_0=\pt_{T(\pt)}$) and for $X, Y \in [S,\Set]$:
\begin{align*}
T_2(X,Y)_s : \sscoprod_{t \in S} \A(t,s) \sstimes  (X\sstimes Y)_t  &\to \sscoprod_{t,t' \in S}  \A(t,s) \sstimes X_t  \sstimes \A(t',s) \sstimes Y_{t'}, \\
(f,x,y) &\mapsto (f,x,f,y).
\end{align*}

A natural transformation 
$\beta\colon \id_\A \sstimes \id_\A \to T \sstimes T$ is given by
$(\beta_{X,Y})_s (x,y) = (a_s, x, b_s, y)$ for a certain choice of $a, b \in \prod_{s \in S} \A(s,s)$.
The corresponding natural transformation $\Phi^\beta_{X,Y}$ is given by:
\begin{align*}
(\Phi^\beta_{X,Y})_s : \sscoprod_{t,u\in S} \A(t,s) \sstimes X_t \sstimes \A(u,t) \sstimes Y_u &\to \sscoprod_{t,u\in S} \A(t,s) \sstimes X_t \sstimes \A(u,s) \sstimes Y_u,\\
(f,x,g,y) &\mapsto (f a_t, x, f b_t g, y).
\end{align*}
The fusion morphism $H^l$ \eqref{eq:8} is obtained when $a_t=\mathrm{id}_t=b_t$ for all $t\in S$.

The category $\A$ is a groupoid if and only if $\A$ is `slack left Hopf':

\begin{thm}\label{thm-cat-monad} Let $\A$ be a small category with set of objects $S$, and let $T$ be the comonoidal monad on $[S,\Set]$ defined above. The following are equivalent:
\begin{enumerate}[(i)]
\item \label{item:grpd1} $T$ is a slack left Hopf monad;
\item \label{item:grpd2} $T$ is a left Hopf monad;
\item \label{item:grpd3} $\A$ is a groupoid.
\end{enumerate}
\end{thm}

\begin{proof}
The equivalence of \ref{item:grpd1} and \ref{item:grpd2} is a special case of Theorem~\ref{th-hopf-cart}.
Now \ref{item:grpd2} means that $H^l$ is an isomorphism, that is:
$$\forall (h : t \to s, k : u \to s), \, \exists! \,\, ( g: u \to t) \blabla{with} h g = k,$$
and this is equivalent to saying that $\A$ is a groupoid (as it implies that any morphism is a split epimorphism), so \ref{item:grpd2} and \ref{item:grpd3} are equivalent.
\end{proof}

\begin{rk}
    The equivalence of \ref{item:grpd2} and \ref{item:grpd3} is shown in a somewhat different form in \cite[Prop.~4.5]{zbMATH06450114}, in the language  of bimonoids in a duoidal category.
\end{rk}

\begin{rk}
Theorem~\ref{thm-cat-monad} asserts here that if $\A$ is a small category such that  there exist  $a, b \in \prod_{s \in S} \A(s,s)$ satisfying
$$\forall (h \colon t \to s, k \colon u \to s), \, \exists! \,\,  (f\colon t \to s, g\colon u \to t) \blabla{with} f a_t =h \blabla{and}  f b_t g = k.$$
then $\A$ is a groupoid, a not very obvious fact which can be proved by hand.
\end{rk}

\subsubsection{Categories as lax magma comonads} There is also an interpretation of a small category as a comonad, and here again the category is a groupoid if and only if it is `slack left Hopf', although that result is not  a consequence of Theorem~\ref{th-hopf-cart}.
Let $\A$ be a small category, with set of objects $S$. As shown by B.~Day and R.~H.~Street in~\cite{Day-Street:quantumcat}, $\A$ induces a monoidal comonad on
$\mathbf{Set}/S\sstimes S$ which is Hopf precisely when the category is a
groupoid, so we may wonder when such a comonad is slack Hopf.

Recall first that if $\C$ is a cartesian category, any object $X$ has a unique structure of a comonoid, with coproduct the diagonal morphism $\Delta_X \colon X \to X \times X$, and therefore defines a comonad
$X \times ?$ on $\C$. Moreover, the forgetful functor $\C/X \to \C$ is comonadic, and its comonad is precisely $X \sstimes ?$.

In particular denote by $A$ the set of arrows of $\A$, viewed as an object of $\Set/S^2$ via $(s,t) : A \to S^2$, and let $C = A\sstimes_{S^2}?$ be the corresponding comonad on $\Set/S^2$.

Now equip $\Set/S^2$ with the monoidal structure given by composition of
spans: if we regard $U$ and $U'$ as directed graphs, then $U\otimes U'$ is the directed graph of pairs of composable edges.
This tensor product has unit
object  $\Delta_S\colon S\to S^2$.

The comonad $C$ associated with $\A$ is monoidal as follows. For $U,U'$ in $\Set/S^2$ the natural
morphism $C(U)\otimes C(U')\to C(U\otimes U')$ has the form
\begin{equation*}
  (A\sstimes_{S^2}U)\sstimes_S(A\sstimes_{S^2}U')
  \to
  A\sstimes_{S^2}(U\sstimes_S U'),
  \label{eq:66}
  \qquad
  (a,u,a',u')\mapsto(a'\cdot a,u,u')
\end{equation*}
where $a'\cdot a$ is the composition of $a$ and $a'$ in the category $\A$.
The morphism $\Delta_S\to C(\Delta_S) = A\sstimes_{S^2}\Delta_S$ is given by $x\mapsto (1_x,x)$.

As monoidal comonads are comonoidal monads on the opposite category, it makes sense to look for slack Hopf structures for $C$. Note that we can not apply Theorem~\ref{th-hopf-cart} in its dual form because, although the category $\Set/S^2$ is cocartesian, the composition of spans is not its coproduct.

A slack (left or right) Hopf structure for the comonad $C$ is given by a natural
transformation
$$\beta_{U,U'}\colon C(U)\otimes C(U')\to U\otimes U'$$
satisfying an extra property (that a certain morphism induced by $\beta$ should
be invertible).
The following lemma shows that there is no wriggle room for such a $\beta$.

\begin{lemma}\label{lemma:rigid}
There is only one natural transformation $C(U) \otimes C(U') \to U \otimes U'$, namely $\eps \otimes \eps$, where $\eps$ is the counit of the comonad $C$.
\end{lemma}

\begin{proof}
Let $u \in U$ and $u' \in U'$, and denote $(x,y)$ (resp. $(x',y')$) the image of $u$ (resp. $u'$) in $S^2$. This defines two objects $\pt_u = (x,y) : \pt \to S^2$ and $\pt_{u'} = (x',y') : \pt \to S^2$ in $\Set/S^2$.
By naturality of $\beta$ applied to the morphism $(u,u') : \pt_u \otimes \pt_{u'} \to U \otimes U'$,
we obtain $\beta_{U,U'} (a,u,a',u') = (u, u')$, that is, $\beta = \eps \otimes \eps$.
\end{proof}

To sum things up:

\begin{thm}\label{thm-cat-comonad} Let $\A$ be a small category with set of objects $S$, and let $C$ be the monoidal comonad on $\Set/S^2$ defined above. The following are equivalent:
\begin{enumerate}[(i)]
\item $C$ is a slack left Hopf comonad;
\item $C$ is a left Hopf comonad;
\item $\A$ is a groupoid.
\end{enumerate}
Moreover if these conditions hold $C$ has exactly one slack Hopf structure.
\end{thm}

\section{Slack Hopf comagma algebras}
\label{sec:examples}


%

\subsection{Slack Hopf comagma monoids}
\label{sec:slack-hopf-comagm}

Let $\B$ be a braided monoidal category, with braiding $\tau$. Given two monoids $(M,m,u)$ and $(M',m',u')$ in $\B$, we endow $M \otimes M'$ with a monoid structure,  with product
$(m \otimes m')(M \otimes \tau_{M',M} \otimes M')$ and unit $u \otimes u'$. This defines a monoidal structure on the category of monoids in $\B$.

A \emph{comagma monoid in $\B$} is a comagma in the category of monoids in $\B$, that is, a monoid $M$ in $\B$ endowed with a monoid morphism $\Delta : M \to  M \otimes M$.

Let $(M,\Delta)$ be a comagma monoid in $\B$. Then $\Delta$ defines a colax magma structure $T_2$ on the monad $M \otimes ? $ on $\B$, defined by
$$T_2(X,Y) = (M \otimes \tau_{X,M} \otimes Y)(\Delta \otimes X \otimes Y).$$

Now let $(M,\Delta)$ be a comagma monoid in $\B$. A morphism $v : \un \to M \otimes M$ defines a morphism $H^v : M \otimes M \to M \otimes M$ as follows:
$$H^v = (m \otimes m)(M \otimes \Delta \otimes m)(M \otimes \tau_{M,M} \otimes M)(v \otimes M \otimes M).$$
We say that $v$ is a \emph{slack left Hopf structure on $M$} if $H^v$ is an isomorphism. If such is the case, $v$ defines a slack left Hopf structure $\beta$ on the colax magma monad $M \otimes ?$, given by:
$$\beta_{X,Y} = (M \otimes \tau_{M,X}\otimes Y)(v \otimes X \otimes Y).$$
Note that the comagma monad $M \otimes ?$ may have slack left Hopf structures that are not of this form.

A \emph{slack left Hopf monoid in $\B$} is a comagma monoid which admits a slack left Hopf structure.

\begin{exa}
Let $M$ be an ordinary monoid. Then the monad $M \times ?$ on $\Set$ has a unique comonoidal structure, which corresponds with the comagma monoid structure defined by the diagonal
$\Delta : M \to M \times M$. By Theorem~\ref{th-hopf-cart}, $M \times ?$ is a slack left Hopf monad if and only if it is a left Hopf monad, and this is equivalent to saying that the fusion map
$H : M^2\to M^2,  (x,y) \mapsto (x,xy)$
is a bijection, equivalently, that $M$ is a group. This can be summarised as `a slack Hopf monoid in $\Set$ is a group'.
A slack left Hopf structure on $M \times ?$ is of the form $(x,y) \mapsto (a,x,b,y)$, with $a, b \in M$. Consequently, the following apparently incongruous  statement holds: if there exists $(a, b) \in M^2$ such that the map $M^2 \to M^2,(x,y) \mapsto (xa,xby)$
is a bijection, then $M$ is a group. Proving it by hand is a little tricky, and is left to the reader as an exercise.
\end{exa}

\begin{rk}
Let $\C$ be a monoidal category, and let $M$ be a comagma monoid in the categorical center  $\Z(\C)$ of $\C$. Then $M$ defines a colax magma monad on $\C$ as follows. One may view $M$ as a monoid in $\C$ with a compatible half-braiding $\sigma$, and the monad $M \otimes ? : \C \to \C$ has a colax magma structure  $T_2(X,Y) = (M \otimes \sigma^{-1}_{X,M} \otimes Y)(\Delta \otimes X \otimes Y)$. Moreover, if $M$ is a slack left Hopf monoid in $\Z(\C)$, then $M \otimes ?$ is a slack left Hopf monad on $\C$.
\end{rk}

\subsection{Comagma Algebras}\label{sec:comagma-algebras}

Let $\kk$ be a commutative ring. The monoidal category of $\kk$-modules $(\Mod_\kk, \otimes = \otimes_\kk, \kk)$ is  equivalent to the category of right exact $\kk$-linear endofunctors of $\Mod_\kk$  via the functor $E \mapsto E \otimes ?$. Thus, a $\kk$-linear right exact slack left Hopf monad on $\Mod_\kk$ can be interpreted as a $\kk$-module endowed with additional structures.

A \emph{comagma algebra over $\kk$} is a comagma monoid in the braided category $\Mod_\kk$ (with its unique braiding). In other words, it is a $\kk$-algebra $A$ endowed with an algebra morphism $\Delta: A \to A\otimes A$.
A comagma algebra structure on a $\kk$-module $A$ is the same thing that a colax magma monad structure on the endofunctor $A \otimes \, ?$.

\begin{rk}\label{rk:sweedler}
For concision, we will use the following variant of Sweedler's notation. If $x$ is an element of a comagma algebra $A$ we write
$\Delta(x) = x_{(1)} \otimes x_{(2)}$. \footnote{We decided to cut ourselves some slack and omit the summation sign.}
Since $\Delta$ is not assumed to be coassociative, we write for instance $$(\Delta \otimes A)\Delta(x) = x_{(1)(1)}\otimes x_{(1)(2)} \otimes x_{(2)}.$$ 
In addition, for an element $u \in A^{\otimes n}$ we adopt the notation $u = u^{(1)} \otimes \dots \otimes u^{(n)}$, and if $u$ is invertible,
$u^{-1} =  u^{(-1)} \otimes \dots \otimes u^{(-n)}$.
\end{rk}

If $A$ is a $\kk\ti$algebra we denote by $A^e$ the enveloping algebra $A^\opp \otimes A$, and by $\cdot$ its product.

\subsection{Slack Hopf Algebras}
\label{sec:slack-Hopf-algebras}

A \emph{slack left Hopf algebra} is a slack left Hopf monoid in $\Mod_\kk$, that is, a comagma algebra $A$ with a slack left Hopf structure $v \in A \otimes A$. The condition on $v$ is that the morphism
$$H^v : A \otimes A \to A \otimes A, \quad x \otimes y \mapsto x_{(1)}v^{(1)}\otimes x_{(2)}v^{(2)}y$$ is an isomorphism (with the conventions of  Remark~\ref{rk:sweedler}).

We denote by $H^l$ the \emph{left fusion morphism} $H^l = H^{1\otimes 1} = (A \otimes m)(\delta \otimes A)$.

Note that a comagma algebra over $\kk$ is a slack left Hopf algebra if and only if the associated colax magma monad is a slack left Hopf monad.

\begin{rk}\label{rk:shbimod}
Let $(A,\Delta)$ be a comagma algebra. Consider two structures of $A-A$ bimodule on $A \otimes A$: $_\bullet A \otimes A_\bullet$, defined by $a\cdot (x \otimes y) \cdot b = ax \otimes yb$, and $_\Delta A \otimes A_\bullet$, defined by $a\cdot (x \otimes y) \cdot b =  \Delta(a)(x \otimes yb)$. Then the assignment $v \mapsto H^v$ defines a bijection between elements of $A \otimes A$ and $A-A$ linear maps $_\bullet A \otimes A_\bullet \to _\Delta\! A \otimes A_\bullet$. As a result, $A$ is slack left Hopf if and only if the $A-A$ bimodules $_\bullet A \otimes A_\bullet$ and $_\Delta A \otimes A_\bullet$
are isomorphic.
\end{rk}

\begin{exa}
A Hopf algebra over $\kk$ is a slack left Hopf algebra with slack left Hopf structure $v=1 \otimes 1$.
\end{exa}
\begin{rk}\label{rk:near-com}Let $A$ be a comagma algebra and assume $$(A \otimes m)(\Delta \otimes A) = (A \otimes m^\opp)(\Delta \otimes A),$$
as is the case in particular when $A$ is commutative or if $\Delta = A \otimes u$. Then for $v \in A\otimes A$, $$H^v = v \cdot H^l,$$
where $\cdot$ denotes the product of $A^e = A^\opp \otimes A$.
\end{rk}

\begin{exa}
Let $(A,\Delta)$ be a comagma algebra, and let $t \in (A \otimes A)^*$. Define $\Delta_t$ to be the conjugate of $\Delta$ by $t$, that is $\Delta_t = t \Delta t^{-1}$.
For $v \in A\otimes A$ denote by $H^v_t$ the morphism $H^v$ relative to the comagma algebra $(A,\Delta_t)$.
One computes easily:
$$H^v_t = t\, H^{t^{-1}v}.$$
In particular, $v$ is a slack left Hopf structure for $(A,\Delta_t)$ if and only if $t^{-1}v$ is one for $(A,\Delta)$.
\end{exa}

\begin{exa}\label{ex:alg-sh}
Let $A$ be a $\kk$-algebra, and let $\Delta = A \otimes u$, where $u$ is the unit of $A$. Then $(A,\Delta)$ is a slack left Hopf algebra, with slack left Hopf element $1 \otimes 1$. (This is a special case of Example~\ref{ex:alg-slack}).
Let $t \in A\otimes A$.
In the notations of the previous example, and using Remark~\ref{rk:near-com}, we have for $v \in A \otimes A$:
$$H^v_t(\xi) = t\,\bigl( (t^{-1}v) \cdot \xi\bigr).$$
As a result, $v$
is a slack left Hopf structure for $(A,\Delta_t)$ if and only if  $t^{-1}v$ is invertible in $A^e = A^\opp \otimes A$. In particular $t$ is a slack left Hopf structure, so  $(A,\Delta_t)$ is a slack left Hopf algebra.
\end{exa}
\begin{exa}\label{exa:mat} Let $\kk$ be a field and let $A$ be the matrix algebra $M_n(\kk)$. All comagma structures on $A$  are of the form $\Delta_t$ in the notation of the previous example, since any two algebra morphisms between two matrix algebras are conjugate to each other.
Let $t\in A \otimes A \iso \End(\kk^n \otimes \kk^n)$ be the invertible element defined by $x \otimes y \to y \otimes x$.
Then $t^{-1}= t$ is not invertible in $A^e$ if $n > 1$ (left as an exercice), and therefore $(M_n(\kk),\Delta_t)$ is a slack left Hopf algebra whose left fusion morphism $H^l$ is not an isomorphism for $n > 1$.
\end{exa}

\begin{lemma}\label{lemma:commut}
Let $A$ be a commutative comagma algebra over a  commutative ring $\kk$.
Then $A$ is slack left Hopf if and only if the left fusion morphism $H^l$ is an isomorphism.
Moreover, in that case the slack left Hopf structures are exactly the invertible elements of $A^e$.
\end{lemma}

\begin{proof} Let $v \in A^e$. Since $A$ is commutative $H^v = v \cdot H^l$ by Remark~\ref{rk:near-com}.
In particular, if $v$ is a slack left Hopf structure, $1 \otimes 1  = H^v(w) = v \cdot H^l(w)$, so $v \in (A^e)^\inv$ because $A^e$ is commutative, and $H^l$ is an isomorphism because so are $H^v$ and $v \,\cdot\, ?$. Moreover if $H^l$ is an isomorphism,
$H^v$ is an iso if and only if $v$ is invertible in $A^e$.
\end{proof}

\subsection{Uniqueness of slack Hopf structures.}
Let $A$ be a slack left Hopf algebra, with corresponding slack left Hopf monad $T = A \otimes\, ?$. The monoid $\M^l(T)$ defined in Section~\ref{sec:unicity} is nothing but the enveloping algebra
$A^e = A^\opp \otimes  A$ of $A$, with its product $\cdot$ given by $(x \otimes y) \cdot (z \otimes t) = zx \otimes yt$. In particular, $\G^l(T) = (A^e)^\inv$. Moreover $\M^l(T)= A^e$ acts on $A \otimes A$ on the right by:
$$v \lhd g = H^v g = g^{(1)}_{\ \ (1)}v^{(1)} \otimes g^{(1)}_{\ \ (2)}v^{(2)}g^{(2)}.$$
If $v \in A\otimes A$ is a slack left  Hopf structure, then $\{v \lhd g; g \in (A^e)^\inv\}$ is the set of all slack left Hopf structures.

\subsection{Internal Homs for modules over a slack Hopf algebra.} \label{sec:closed-alg}
If $A$ is a slack left Hopf algebra, with slack left  Hopf structure $v$, then the category $A-\Mod$ of left $A$-modules is left closed. Using Theorem~\ref{thm-antipode}, we now make the left internal Homs of $A-\Mod$  explicit in terms of $v$.

We introduce new notation. If $V, W$ are two $A$-modules, then $\Hom_\kk$ is a right  $A^e$-module, where $A^e$ is the enveloping algebra $A^\opp \otimes A$, with right action defined, for a linear map $f \in \Hom_\kk(V,W)$ and an element $\xi$ of $A^e$, by
$$(f \lhd \xi)(x) = \xi^{(1)}f(\xi^{(2)} x).$$
Notice that $A^e$ has a faithful right module of the above form. For example
the Hom between the free left $A$-modules $A$ and $A\otimes A$: $((u\otimes
\mathrm{id}_A)\lhd \xi)(1)=\xi$, for $\xi\in A^e$.

Now define $\nabla \colon A \to A \otimes A$ by
$$\nabla(x) = {H^v}^{-1}(v^{(1)}x \otimes v^{(2)}).$$
The slack left Hopf antipode $S^l_{Y,X}: A \otimes \Hom_\kk(Y,A\otimes X) \to \Hom_\kk(Y, A \otimes X)$ for $X, Y$ in $\Mod_\kk$ is defined by $S^l_{X,Y}(a \otimes f)(y) = \nabla^{(1)}(a) \otimes f(\nabla^{(2)}(a) \otimes y)$, so that, for $A$-modules  $V$ and $W$, the action of $A$ on $\lhom{V}{W} = \Hom_\kk(V,W)$ is given by
$a f  = f \lhd \nabla(a)$.

The unit and counit of the internal Hom adjunction are given by:
\begin{align*}
h^{V}_{W} :\,\, & W \to \Hom_\kk(V,W\otimes V), & y \, \mapsto \,&(x \mapsto v^{(1)} \cdot y \otimes v^{(2)}\cdot x),\\
e^{V}_{W} :\,\, & \Hom_\kk(V,W) \otimes V \to W, & f \otimes x \, \mapsto\, &w^{(1)} \cdot f(w^{(2)}\cdot x),
\end{align*}
where $w = {H^v}^{-1}(1 \otimes 1)$.

\begin{lemma}
    \label{l:nablaalg}
    The morphism $\nabla$ is a an algebra antimorphism $A\to A^e$.
\end{lemma}
\begin{proof}
    As pointed out above, $A^e$ has a faithful module of the form $\Hom_{\kk}(V,W)$. It is clear that that $\nabla(ab)$ and $\nabla(b)\nabla(a)\in A^e$ act in the same way on Homs, so they are equal. Likewise for $\nabla(1)$ and $1\otimes 1$.
\end{proof}

\begin{rk}
Two different slack left Hopf structures $v$, $v'= v \lhd \gamma$ (with $\gamma \in (A^e)^\inv$) define left internal Homs $\lhom{V}{W}_v$ and $\lhom{V}{W}_{v'}$, with the same underlying $\kk$-module $\Hom_\kk(V,W)$ but different $A$-module structures, evaluation and coevaluation morphisms. However, the uniqueness of left internal Homs up to unique isomorphism yields an isomorphism of $A$-modules $\lhom{V}{W}_v \iso \lhom{V}{W}_{v'}$, which is given by $f \mapsto f \lhd \gamma$. 
\end{rk}

\subsection{The relation between $(\Delta,v)$ and $(\nabla,w)$}

\begin{rk}\label{rk:adj-eq}
Let $A$ be a slack left Hopf algebra, with slack left Hopf structure $v$. Then, with the above notations, we verify easily:
$${H^v}^{-1}(x \otimes y) = \nabla(x) \cdot w  \cdot (1 \otimes y) =  w^{(1)}\nabla^{(1)}(x)  \otimes \nabla^{(2)}(x) w^{(2)} y.$$
The $A$-linearity of the adjunction unit and counit $h$ and $e$ translates (respectively) as
$$(1) \quad (x \otimes 1) \cdot v  = H^v \nabla(x) \blabla{and} (2)\quad w \cdot (x \otimes 1)  = {H^v}^{-1} \Delta(x) \blabla{(for $x \in A$)}$$

\noindent
and the adjunction identities translate as
$$(3) \quad H^v w = 1 \otimes 1 \blabla{and} (4)\quad {H^v}^{-1}v = 1 \otimes 1.$$
Note that these equations are nothing new: (4) is obvious; (1) and (3) hold by definition of $\nabla$ and $w$; (2) follows from left $A$-linearity of $H^v$.
\end{rk}

\begin{rk}
Given an algebra $A$, there is a bijective correspondance between \begin{enumerate}[(i)]
\item data $(\Delta,v)$, where $\Delta$ is an algebra morphism $A \to A \otimes A$ and $v$ is a slack left Hopf structure, that is, an element of $A \otimes A$ such that
$$H^v = H^{\Delta,v} : A \otimes A \to A \otimes A, \quad x \otimes y \mapsto  \Delta^{(1)}(x) v^{(1)} \otimes \Delta^{(2)}(x) v^{(2)} y$$
is an isomorphism,
\item data $(\nabla,w)$, where $\nabla$ is an algebra morphism
$A^\opp \to A^e$ and $w$ is an element of $A^e$ such that $$Q^{\nabla,w} : A \otimes A \to A \otimes A, \quad x \otimes y \mapsto  w^{(1)} \nabla^{(1)}(x) \otimes \nabla^{(2)}(x)w^{(2)}y$$
is a isomorphism,
\end{enumerate}
which is characterised by the fact that $H^{\Delta,v}$ and $Q^{\nabla,w}$ are mutually inverse.

The reason for this duality is that a left closed magma structure on $\C = A\ti\Mod$ consists in a tensor product functor $\otimes: \C \times \C \to \C$ and a left internal Hom functor $\lhom{-}{-}: \C^\opp \times \C \to \C$, related by the adjunction $V \otimes \, ? \adj \lhom{V}{?}$. We have considered the tensor product as the primary datum, which lead us to encode the left closed magma structure as $(\Delta,v)$. Had we given precedence to the left functor, we would have encoded it as $(\nabla,w)$.
\end{rk}

The following lemma exploits the fact that  $(\Delta,v)$, $(\nabla,w)$ play a symmetric r\^ole. Given a comagma algebra $A$ with slack left Hopf structure $v$, and $s \in A^e$, define
$$Q^s: A \otimes A \to A \otimes A, \quad x \otimes y \mapsto\nabla(x)\cdot s \cdot (1 \otimes y) = s^{(1)} \nabla^{(1)}(x) \otimes \nabla^{(2)}(x)s^{(2)}y.$$

\begin{lemma}\label{lemma:vw-duality}
Let $A$ be a comagma algebra with slack left Hopf structure  $v$. Then:
\begin{enumerate}
\item $v$ is the only element $t$ in $A^e$ satisfying $$H^t \nabla(x) = (x \otimes 1)\cdot t \blabla{and} H^t(w) = 1 \otimes 1;$$
\item $w$ is the only element $s$ in $A^e$ satisfying $$Q^s \Delta(x) = s \cdot (x \otimes 1) \blabla{and} Q^s(v) = 1 \otimes 1.$$
\end{enumerate}
\end{lemma}

\begin{proof} Let us prove Assertion (2). That $s = w$ satisfies the stated properties results from Remark~\ref{rk:adj-eq}, observing first that $Q^w = H^{v -1}$. So assume
$Q^s \Delta(x) = s \cdot (x \otimes 1)$ and $Q^s(v) = 1 \otimes 1$. We have then
\begin{align*}
w &= Q^s(v)\cdot w = \nabla(v^{(1)}) \cdot s \cdot (1 \otimes v^{(2)}) \cdot w = \nabla(v^{(1)}) \cdot s \cdot (w^{(1)} \otimes 1) \cdot (1 \otimes v^{(2)}w^{(2)})\\
&=   \nabla(v^{(1)})
 \cdot \nabla(w^{(1)}_{\,\,(1)}) \cdot s \cdot (1 \otimes w^{(1)}_{\,\,(2)}) \cdot (1 \otimes  v^{(2)}w^{(2)}) \\
& = \nabla(w^{(1)}_{\,\,(1)} v^{(1)}) \cdot s \cdot (1 \otimes w^{(1)}_{\,\,(2)} v^{(2)}w^{(2)})\\ &= \nabla(H^v(w)^{(1)}) \cdot s \cdot (1 \otimes H^v(w)^{(2)})
 = Q^s H^v(w) = Q^s(1 \otimes 1) = s.
\end{align*}
The proof of Assertion (1) goes along the same lines. We know that $t = v$ satisfies the stated properties, and if $t$ satisfies them, a similar computation shows  $v = H^t(w) \,v$ (product in $A \otimes A$) $= H^t Q^w(v) =H^t(1 \otimes 1) = t$.
\end{proof}

\subsection{Slack Hopf algebras with counits}

A \emph{left counit} for a comagma algebra $A$ is an algebra morphism $\eps: A \to \kk$ such that $(\eps \otimes A)\Delta = \id_A$.

\begin{lemma}\label{lemma-comag-unit} Let $A$ be a slack left Hopf algebra with slack left Hopf structure $v$ and a left counit $\eps$.
With our usual notation set $\sigma = (\eps \otimes A)\nabla : A \to A$, $\alf = (\eps\otimes A)w$ and $\bee = (\eps \otimes A)v$. Then $\sigma : A \to A$ is an antimorphism of algebras, and we have:
\begin{enumerate}
\item $\sigma(x_{(1)}) \,\alf \, x_{(2)} = \alf\, \eps(x)$ for $x \in A$;
\item $\nabla^{(1)}(x) \,\bee\, \nabla^{(2)}(x) = \bee\, \eps(x)$ for $x \in A$;
\item $w^{(1)} \,\bee\, w^{(2)} = 1$;
\item $\sigma(v^{(1)}) \,\alf \, v^{(2)} = 1$.
\end{enumerate}
\end{lemma}

\begin{proof} Results immediately from Lemma~\ref{l:nablaalg} and Remark~\ref{rk:adj-eq}: $\sigma$ is an antimorphism of algebras because  $\nabla : A \to A^e$ and $\eps \otimes A : A^e \to A$ are respectively an algebra antimorphism and an algebra morphism. Besides, Equations  (1)--(4) of the Remark, postcomposed by $(\eps \otimes A)$, yield Assertions (1) -- (4) of the Lemma.
\end{proof}

\subsection{Slack Hopf bialgebras and Hopf algebras}

We now consider the case of bialgebras. Recall that if $C$ is a coalgebra and $B$ is an algebra, $\Hom_\kk(C,B)$ is an algebra with convolution product
$*$ defined by $f*g = m_B (f \otimes g) \Delta_C$ and unit $u_B \eps_C$.

Let $A$ be a bialgebra with counit $\eps$. An antipode $S$ for $A$ is an inverse of $\id_A \in \Hom_\kk(A,A)$ for the convolution product. The existence of an antipode is equivalent to the invertibility of the left fusion operator $H^l = (A \otimes m)(\delta \otimes A): A \otimes A \to A \otimes A$. A \emph{Hopf algebra} is a bialgebra with an antipode.

The commutative case is particularly simple,  being a special case of
Lemma~\ref{lemma:commut}:

\begin{cor}\label{cor:commut bialg}
A commutative bialgebra is a Hopf algebra if and only if it is slack left Hopf.\qed
\end{cor}

We now deal with the general case. If $\kk$ is a field,
a $\kk$-algebra $B$ is \emph{profinite} if it is cofiltered limit of finite dimensional algebras. If $B$ is a profinite algebra and $C$ is a coalgebra, then $\Hom_\kk(C,B)$ is a profinite algebra for the convolution product, because a $\kk$-coalgebra is the filtered union of its finite dimensional subcoalgebras.



\begin{thm}\label{thm-slack-Hopf-bialg} Let $A$ be a bialgebra over a commutative ring $\kk$, let $v \in A \otimes A$ be a slack left Hopf structure, and, in the notations of Section~\ref{sec:closed-alg}, set $\alf = (\eps \otimes A)w \in A$ and $\sigma = (\eps \otimes A)\nabla$.
Then:
\begin{enumerate}
\item If $A$ is a Hopf algebra, $\alf$ is invertible;
\item If $\alf$ is invertible, the linear map  $S : A \to A$  defined by
$$x \mapsto  S(x) = \alf^{-1} \cdot \sigma(x)\cdot \alf$$
is a left inverse of $\id_A$  in the convolution algebra $\Hom_\kk(A,A)$;
\item If $\alf$ is invertible and in addition $\kk$ is a field and $A$ is profinite,
$S$ is an antipode and so $A$ is a Hopf algebra.
\end{enumerate}
\end{thm}
\begin{proof}
The bialgebra $A$ is a  Hopf algebra if and only if  $H^l$ is an isomorphism, which is the case if and only if $w$ is invertible in $A^e$ by Corollary~\ref{cor-slh-h}. Since $\eps$ is an algebra morphism, we have that
$\alf = (\eps \otimes A)w$ is invertible in $A$, which proves Assertion~(1).
Assume $\alf$ is invertible and set
$S = \alf^{-1}\sigma\alf$. By Lemma~\ref{lemma-comag-unit}, $\sigma(x_{(1)})\alf x_{(2)} = \alf \eps(x)$.  Then $S(x_{(1)}) x_{(2)} = \eps(x)$, which means that $S$ is left convolution inverse to $\id_A$,
Now, if $\kk$ is a field and $A$ is profinite, so is the convolution algebra $\Hom_\kk(A,A)$. Assertion~(3) results from the fact that in a profinite algebra an element which has an inverse on one side is invertible, and therefore $S$ is an inverse of $\id_A$, that is, an antipode. Note that $(\sigma,\alf,\alf^{-1})$ is a quasi-antipode, see Section~\ref{sec:qha} below.
\end{proof}


\section{Quasi-Hopf algebras}\label{sec:qha} Let $(A,\Delta,\eps,\phi)$ be a quasi-bialgebra over a field $\kk$.
Recall that a \emph{quasi-antipode} is a triple $(\Es,\alf,\bee)$, where $\Es : A \to A$ is an antimorphism of algebras and  $\alf,\bee$ are elements of $A$ satisfying the following four axioms:
\begin{enumerate}[({QA}1)]
\item $\Es(x_{(1)}) \,\alf\, x_{(2)} = \alf\, \eps(x)$, $x \in A$;
\item $x_{(1)} \,\bee\, \Es(x_{(2)}) = \bee\, \eps(x)$, $x \in A$;
\item $\phi^{(1)} \,\bee\, \Es(\phi^{(2)}) \,\alf\, \phi^{(3)} = 1$;
\item $\Es(\phi^{(-1)})  \,\alf\, \phi^{(-2)} \,\bee\, \Es(\phi^{(-3)}) = 1$,
\end{enumerate}
where $\phi = \phi^{(1)} \otimes \phi^{(2)} \otimes \phi^{(3)}$ and $\phi^{-1} = \phi^{(-1)} \otimes \phi^{(-2)} \otimes \phi^{(-3)}$, in agreement with our conventions.

A \emph{quasi-Hopf algebra} is a quasi-bialgebra which admits a quasi-antipode.

Let $A$ be a quasi-Hopf algebra over $\kk$, with quasi-antipode $(\Es,\alf,\bee)$.
Following Schauenburg~\cite{Schau:CFQHA}, define $v, w\in A\otimes A$ by:
    $$
    v=\phi^{(-1)}\otimes\phi^{(-2)}\,\bee\, \Es( \phi^{(-3)}), \quad w =\phi^{(1)} \otimes \Es(\phi^{(2)})\,\alf\,\phi^{(3)}.
    $$
\label{sec:quasi-hopf-slack}
Schauenburg shows that the morphism we denote by $H^v$ is bijective, with
$$H^{v -1}(x \otimes y) = w^{(1)} x_{(1)} \otimes \Es(x_{(2)}) w^{(2)} y$$
 (in particular, $w = H^{v -1}(1 \otimes 1)$, in agreement with our notation).
In other words, $v$ is a slack left Hopf structure, and a quasi-Hopf algebra is slack left Hopf.


\begin{rk}
A slack (right) Hopf structure for a quasi-Hopf algebra is already present in Drinfeld~\cite{QHA}, Proposition 1.5., in a slightly different form since it involves the inverse of $\Es$. \footnote{For a quasi-antipode Drinfeld assumes $\Es$ to be bijective.}
\end{rk}

A slack left Hopf structure $v$ which can be so expressed in terms of a quasi-antipode of $A$ is called a \emph{left Hopf structure}.

\begin{thm}\label{thm-quanti-from-slh} If $v$ is the left Hopf structure associated with a quasi-antipode $(\Es,\alf,\bee)$ then:
\begin{enumerate}
\item $\alf = (\eps \otimes A)w = (\eps \otimes A) {H^v}^{-1}(1 \otimes 1)$;
\item $\bee = (\eps \otimes A)v$;
\item $\Es= \sigma$;
\item $\nabla = (A \otimes \Es)\Delta = (A\otimes \sigma)\Delta$.
\end{enumerate}
\end{thm}
\begin{proof}
That Schauenburg's $w$ is equal to  ${H^v}^{-1}(1 \otimes 1)$ results from the expression given for ${H^v}^{-1}$ and the first two assertions result from the fact that $(\eps \otimes A \otimes A) \phi^{-1}=1 \otimes 1 = (\eps \otimes A \otimes A) \phi$ and $\Es$ is an algebra antimorphism. In order to verify the identity of Assertion (4), where $\nabla(x) = {H^v}^{-1} (v^{(1)} x \otimes v^{(2)})$,
we compute
\begin{multline*}
H^v(A \otimes \Es)\Delta(x) = H^v(x_{(1)} \otimes \Es(x_{(2)})= x_{(1)(1)}\phi^{(-1)}\otimes x_{(1)(2)}\phi^{(-2)} \bee \Es(\phi^{(-3)}) \Es(x_{(2)})\\
=x_{(1)(1)}\phi^{(-1)}\otimes x_{(1)(2)} \phi^{(-2)} \bee \Es(x_{(2)}\phi^{(-3)})
=\phi^{(-1)}x_{(1)} \otimes\phi^{(-2)}x_{(2)(1)}  \bee \Es(\phi^{(-3)} x_{(2)(2)})\\
=\phi^{(-1)}x_{(1)} \otimes\phi^{(-2)}x_{(2)(1)}  \bee \Es(x_{(2)(2)})\Es(\phi^{(-3)}) =\phi^{(-1)}x_{(1)} \otimes \phi^{(-2)}\bee \eps(x_{(2)})\Es(\phi^{(-3)})\\
=\phi^{(-1)}x \otimes\phi^{(-2)}\bee \Es(\phi^{(-3)}) = v^{(1)}x \otimes v^{(2)},
\end{multline*}
from which we conclude that $\nabla = (A \otimes \Es)\Delta$, so $\sigma = (\eps \otimes A) \nabla = \Es$, that is Assertion (3), and so  Assertion (4) holds as well.
\end{proof}

\begin{rk}In particular, the left internal Homs in $A\ti\Mod$ constructed in \cite{NonAss} from a quasi-antipode coincide with the ones defined by the corresponding slack left Hopf structure $v$ (with same module structure and adjunction units and counits). Note also that the axioms of a quasi-antipode match exactly the properties of slack left Hopf structures listed in Lemma~\ref{lemma-comag-unit}.
\end{rk}

How could we characterise left Hopf structures among  slack left Hopf structures? The next section is devoted to this question.


\subsection{Slackness and comparison morphisms}

Given a slack left Hopf structure $v$ for a quasibialgebra $A$, recall that $w = H^{v -1}(1 \otimes 1)$, $\nabla(x) = {H^v}^{-1} (v^{(1)} x \otimes v^{(2)})$ and
$\sigma = (\eps \otimes A)\nabla$. Set $\alf = (\eps \otimes A)w$, $\overline{w} = \phi^{(1)} \otimes \sigma(\phi^{(2)}) \,\alf\, \phi^{(3)} \in A \otimes A$, and define
$$\SL(v)  = \overline{w}^{(1)} v^{(1)}_{\,\,(1)} \otimes \sigma (v^{(1)}_{\,\,(2)}) \overline{w}^{(2)} v^{(2)} \in A^e.$$
The element $\SL(v)$ is called the \emph{slackness} of $v$.

Let $\C$ be a monoidal left closed category, let $X,Y$ be objects of $\C$ and consider the canonical morphism, called \emph{comparison morphism}:
$$c_{X,Y} : Y \otimes \lhom{X}{\un} \to \lhom{X}{Y}$$
obtained by universal property of the internal Hom from the morphism
$$
(Y \otimes \lhom{X}{\un}) \otimes X \xrightarrow{\cong}
Y \otimes (\lhom{X}{\un} \otimes X) \xrightarrow{\quad Y \otimes e^X_\un}
Y \otimes \un \cong Y,
$$

For $E$, $F$ in $\Vect_\kk$ the comparison morphism $F \otimes E^* \to \Hom_\kk(F,E)$ is the canonical map $x \otimes f \mapsto (x \circ f : y \mapsto f(y)x)$ (viewing $x$ as a map $\kk \to E$).

\begin{lemma}\label{lemma-c-d} Given a quasi-bialgebra $A$ with slack left Hopf structure $v$, for any pair of $A$-modules $V, W$ the comparison morphism
$$c_{V,W} : W \otimes \lhom{V}{\un} \to \lhom{V}{W}$$
is given by
$$c_{V,W}(x \otimes f) = (x \circ f)\lhd \SL(v).$$
\end{lemma}

\begin{proof} Taking into account the associativity constraint $a$ of $A\ti\Mod$, the morphism $c_{V,W}$ is obtained by universal property of the internal Hom from the morphism
$$
(W \otimes \lhom{V}{\un}) \otimes V \xrightarrow{a_{W,\lhom{V}{\un},V}}
W \otimes (\lhom{V}{\un} \otimes V) \xrightarrow{W \otimes e^V_\un}
W \otimes \un \cong W,
$$
that is, $c_{V,W} =  \lhom{V}{(W \otimes e^V_\un)a_{W,\lhom{V}{\un},V}} h^V_{W \otimes \lhom{V}{\un}}$.
Using the explicit formulae for internal Homs in $A-\Mod$ and the corresponding adjunction evaluation and coevaluation, as well as the expression of the associativity constraint $a$ in terms of $\phi$, we obtain by straightforward computation firstly
$$(W \otimes e^V_\un)a_{W,\lhom{V}{\un},V} (x \otimes f \otimes y) = \overline{w}^{(1)}x \otimes f(\overline{w}^{(2)}y),$$
and secondly
$$c_{V,W}(x \otimes f)(y) = \overline{w}^{(1)}v^{(1)}_{\,\,(1)} x \otimes f(\sigma(v^{(1)}_{\,\,(2)}) \overline{w}^{(2)}v^{(2)} y)
= \SL(v)^{(1)} x \otimes f(\SL(v)^{(2)} y)$$
as announced.
\end{proof}

\begin{lemma}\label{lemma:sigma-delta} For a slack left Hopf algebra $A$ with slack left Hopf structure $v$, we have
$$sl(v) \cdot \Delta ) ((A \otimes \sigma)\Delta)\cdot \SL(v).$$
\end{lemma}

\begin{proof}
This statement results from Lemma~\ref{lemma-c-d} an expresses the $A$-linearity of $c_{V,W}: W \otimes \lhom{V}{\un} \to \lhom{W}{V}$, using the fact that $(\eps \otimes \sigma)\Delta$ and $\nabla$ encode the action of $A$ respectively on the source and the target. For this we use the fact that one may choose $V$ and $W$ so that
$W \otimes \lhom{V}{\un}$ is a faithful right $A^e$-module, which is the case if $\kk$ is a field, taking $V = W = A$.
\end{proof}

\begin{rk}\label{rk-sl_G0}
For a slack left Hopf structure $v$, the comparison map $c_{\un,W} : W^* \to W^*$ is the identity, from which results:
$(\eps \otimes A)\,\SL(v) = 1$.
\end{rk}

\begin{lemma}\label{lemma-wbar}\label{lemma-KH} Let $A$ be a quasi-biagebra with slack left Hopf structure $v$. For $x\in A$, $\overline{w} \cdot(x \otimes 1) = \overline{w}^{(1)}x_{(1)(1)}\otimes \sigma(x_{(1)(2)})\overline{w}^{(2)}x_{(2)}$.
\end{lemma}

\begin{proof}
This identity expresses the  $A$-linearity of the map $(W \otimes e^V_\un)a_{W,\lhom{V}{\un},V}: (W \otimes \lhom{V}{\un}) \otimes V \to W$, $x \otimes f \otimes y \mapsto \overline w^{(1)} \otimes f(\overline w^{(2)} y)$ met in the course of  the proof of Lemma~\ref{lemma-c-d} (and uses the hypothesis that $\kk$ is a field).
\end{proof}

The following lemma clarifies the dependence of $\SL(v)$ on $v$.

\begin{lemma}\label{lem-v-change} If  $\gamma \in (A^e)^\inv$, setting $\gamma_0 = (\eps \otimes A)\gamma$ we have
$$\SL(v\lhd \gamma) = (1 \otimes \gamma_0^{-1})\cdot \SL(v)\cdot \gamma\,.$$
\end{lemma}

\begin{proof}
This results from the interpretation of $\SL(v)$ from Lemma~\ref{lemma-c-d} using Remark~\ref{sec:closed-alg}, or can be verified by hand using  Lemma~\ref{lemma-wbar}.
\end{proof}

\subsection{Slack Hopf versus Quasi-Hopf}

\begin{thm}\label{thm-crit-quanti} Let $A$ be a quasi-bialgebra and let $v$ be a slack left Hopf structure. Let $\alf = (\eps \otimes A)w$, $\bee = (\eps\otimes A)v$ and $\overline{v} = \phi^{(-1)} \otimes \phi^{(-2)} \,\bee\, \sigma(\phi^{(-3)})$.

The following conditions are equivalent:
\begin{enumerate}[(i)]
\item $v$ is a left Hopf structure;
\item $\SL(v) = 1 \otimes 1$;
\item $\nabla = (A \otimes \sigma)  \Delta$ and $w = \overline{w}$;
\item $\nabla = (A \otimes \sigma)  \Delta$ and $v= \overline{v}$.
\end{enumerate}
Moreover if these hold, the quasi-antipode associated to  $v$ is $(\sigma,\alf,\bee)$.
\end{thm}

\begin{proof}
For $t \in A^e$ define $$K^t : A \otimes A \to A \otimes A,  \quad x \otimes y \mapsto  t^{(1)}x_{(1)} \otimes \sigma(x_{(2)}) t^{(2)} y.$$
By definition, $\SL(v) = K^{\overline{w}}(v)$. Note also that if $\nabla = (A \otimes \sigma)\Delta$, $K^t = Q^t$ as defined just before Lemma~\ref{lemma:vw-duality}.

So if (iii) holds, that is,  $\nabla = (A \otimes \sigma)  \Delta$ and $w = \overline{w}$, we have $K^{\overline w} = K^w = Q^w =  H^{v -1}$, by Remark~\ref{rk:adj-eq}.
Then $\SL(v) = H^{v-1}(v) = 1 \otimes 1$ so (iii) $\Rightarrow$ (ii).

If $v$ is a left Hopf structure, by Theorem~\ref{thm-quanti-from-slh} the corresponding quasi-antipode is  $(\sigma,\alf, \bee)$, $v = \overline v$, $w = \overline w$, and $\nabla = (A \otimes \sigma)\Delta$; in particular (iii) holds and so  $\SL(v) = H^{v -1} (v) = 1 \otimes 1$. In other words (i) $\Rightarrow$ (ii), (iii) and (iv). In passing we have proved the `moreover' assertion of the Theorem.

We know from Lemma~\ref{lemma-comag-unit} that, for any slack left Hopf structure, $\sigma$ is an antimorphism of algebras and the following identities hold:
\begin{enumerate}
\item $\sigma(x_{(1)}) \,\alf \, x_{(2)} = \alf\, \eps(x)$, $x \in A$;
\item $\nabla^{(1)}(x) \,\bee\, \nabla^{(2)}(x) = \bee\, \eps(x)$, $x \in A$;
\item $w^{(1)} \,\bee\, w^{(2)} = 1$;
\item $\sigma(v^{(1)}) \,\alf \, v^{(2)} = 1$.
\end{enumerate}
Comparing these with the axioms of a quasi-antipode:
\begin{enumerate}[({QA}1)]
\item $\sigma(x_{(1)}) \, \alf x_{(2)} = \alf\, \eps(x)$, $x \in A$;
\item $x_{(1)} \,\bee \sigma(x_{(2)}) = \bee\, \eps(x)$, $x \in A$;
\item $\phi^{(1)} \,\bee\, \sigma(\phi^{(2)}) \,\alf\, \phi^{(3)} = 1$;
\item $\sigma(\phi^{(-1))})  \,\alf\, \phi^{(-2)} \,\bee\, \sigma(\phi^{(-3)}) = 1$;
\end{enumerate}
we observe that QA1 is for free, and that the conditions $\nabla = (\eps \otimes \sigma)\Delta$,  $w = \overline{w}$  and $v = \overline v$ imply respectively QA2, QA3 and QA4.

Therefore if both (iii) and (iv) hold, $(\sigma, \alf, \bee)$ is a quasi-antipode and $v = \overline{v}$ is the corresponding slack left Hopf structure. In other words the conjunction of (iii) and (iv) implies (i).  At this stage the equivalence of (iii) and (iv) will come in handy.

\begin{lemma}\label{lemma:w-v} If $\nabla = (A \otimes \sigma) \Delta$, then $w = \overline w$ if and only if $v = \overline v$.
\end{lemma}

\begin{proof}
We use Lemma~\ref{lemma:vw-duality}, and with our assumption on $\nabla$, we have $Q^s = K^s$ and QA1, QA2 hold. Now if $w = \overline w$, QA3 holds.
This implies $H^{\overline v}(\overline w) = 1 \otimes 1$. Indeed,
\begin{multline*}
H^{\overline v}(\overline w)= (A \otimes (m^5(A \otimes \bee \otimes \sigma \otimes \alf \otimes A)) ( (\phi^{-1} \otimes 1)(\Delta \otimes A \otimes A)\phi)\\
\hfill
= (A \otimes m^5(A \otimes \bee \otimes \sigma \otimes \alf \otimes A))((A \otimes A \otimes \Delta)\phi^{-1} (1 \otimes \phi)(A \otimes \Delta \otimes A)\phi)\hfill \text{(5-gon)}\\
\hfill
= 1 \otimes \phi^{(1)} \bee \phi^{(2)}  \alf \phi^{(3)}
\hfill \text{(using QA1 and QA2)} \\
\hfill = 1 \otimes 1  \hfill \text{(by QA3)}.
\end{multline*}
So $H^{\overline v}(w) = H^{\overline v}(\overline w) = 1 \otimes 1$. By Remark~\ref{rk:adj-eq}, Equation (1), we have $\overline w \cdot (x \otimes 1)
= {H^v}^{-1} \Delta(x)$, and applying Assertion (1)
of Lemma~\ref{lemma:vw-duality}, we get $v = \overline v$.

The converse is proved in a similar way: if $v = \overline v$, QA4 holds and this implies $K^{\overline w}\overline v = 1 \otimes\sigma(\phi^{(-1))})  \,\alf\, \phi^{(-2)} \,\bee\, \sigma(\phi^{(-3)}) = 1 \otimes  1$.  Using Equation (2) of Remark~\ref{rk:adj-eq} and  Assertion (2) of Lemma~\ref{lemma:vw-duality}, we get $w = \overline w$.
\end{proof}

Therefore (iii) $\Leftrightarrow$ (iv) and at this point it is enough to show (ii) $\Rightarrow$ (iii),
%
%
%
so assume  $\SL(v) = 1 \otimes 1$.
Lemma~\ref{lemma:sigma-delta} yields  $\nabla = (A \otimes \sigma) \Delta$, and so QA2 holds. We have  $Q^{\overline w}(v) = K^{\overline w}(v)  =\SL(v) = 1 \otimes 1$. By Lemma~\ref{lemma-wbar} we also have $\overline w \cdot (x \otimes 1) =  K^{\overline w}\Delta(x) = Q^{\overline w} \Delta(x)$.
Therefore by Assertion (2) of Lemma~\ref{lemma:vw-duality}, $w = \overline{w}$. In other words (ii) $\Rightarrow$ (iii), which concludes the proof.
%
\end{proof}

\begin{cor}\label{cor-quanti}
Let $A$ be a slack left Hopf quasi-bialgebra. The following are equivalent:
\begin{enumerate}[(i)]
\item $A$ is a quasi-Hopf algebra;
\item There exists a left Hopf structure;
\item There exists a slack left Hopf structure $v$ such that $\SL(v)$ is invertible in $A^e$;
\item For any slack left Hopf structure $v$, $\SL(v)$ is invertible in $A^e$.
\end{enumerate}
\end{cor}

\begin{proof}
By definition we have (i) $\Leftrightarrow$ (ii). By Theorem~\ref{thm-crit-quanti}, a left Hopf structure is nothing but  a slack left Hopf structure $v$ such that $\SL(v) = 1\otimes 1$, so (ii) $\Rightarrow$ (iii).  Now let $v$, $v'$ be two slack left Hopf structures. There exists $\gamma \in (A^e)^\inv$ unique such that
$v' = v \lhd \gamma$, and by Lemma~\ref{lem-v-change}, $\SL(v') = (1 \otimes \gamma^{-1}_0)\cdot \SL(v) \cdot \gamma$, where $\gamma_0 = (\eps \otimes A)\gamma \in A^\inv$. In particular if $\SL(v)$ is invertible in $A^e$ so is $\SL(v')$, so (iii) $\Leftrightarrow$ (iv).
Moreover if $\SL(v)$ is invertible, let $\gamma$ be its inverse in $A^e$. We have $\gamma_0 = 1$ by Remark~\ref{rk-sl_G0}, so $\SL(v') = 1 \otimes 1$, which means $v$ is a left Hopf structure hence (iii) $\Rightarrow$ (ii).
\end{proof}

\begin{rk}
The main difficulty of Theorem~\ref{thm-crit-quanti} -  Corollary~\ref{cor-quanti} is to prove that if a quasi-bialgebra admits a slack left Hopf structure $v$ satisfying $\SL(v) \in (A^e)^\inv$, it is quasi-Hopf.  There is an alternative proof of this fact when $A$ is profinite,  based on classical results about duals of finite dimensional modules on a quasi-bialgebra, which goes as follows. Firsty, it is well-known that a profinite quasi-bialgebra $A$ is quasi-Hopf if and only if the category $A\ti\mod$ of finite dimensional $A$-modules is left autonomous, and the corresponding forgetful functor
$U_0 : A\ti\mod \to \vect_\kk$ is left autonomous, that is, endowed with an isomorphism of functors $U_0(\ldual{V}) \iso U_0(V)^*$. Secondly, it is also well-known that a monoidal left closed category $\C$ is
left autonomous if and only if for any object $V$ the comparison morphism $c_{V,V}$ is an isomorphism, and in that case $\ldual{V} = \lhom{V}{\un}$. Lastly, we know from Lemma~\ref{lemma-c-d} that this is the case if $\SL(v)$ is invertible, and then $U_0$ is left autonomous because it is slack left closed.
\end{rk}

\subsection{A categorical characterisation of quasi-Hopf algebras}

Let $U : \D \to \C$ be a colax magma functor between monoidal left closed categories equipped with a morphism $U_0: U(\un) \to \un$. Let $\mathfrak{E}$ be a slack left closed structure on $U$, that is, a collection of isomorphisms $$\mathfrak{E}^A_B: U\lhom{A}{B} \iso \lhom{UA}{UB}$$
natural in $A$, $B$ in $\D$.

Then  $(U,\mathfrak{E})$ is said to preserve comparison morphisms if for any pair $A, B$ of objects of $\D$ the following diagram commutes:
$$
\begin{tikzcd}
    U(B) \otimes U\lhom{A}{\un}
    \ar[d,"UB \otimes\mathfrak{E}^A_\un "']
    &[30pt] U(B \otimes \lhom{A}{\un}) \ar[r,"U(c_{A,B})"]\ar[l,"{U_2(B,\lhom{A}{\un})}"']
    &[25pt] U\lhom{A}{B}
    \ar[d,"\mathfrak{E}^A_B"]
    \\
    UB \otimes \lhom{UA}{U\un}
    \ar[r,"UB \otimes \lhom{UA}{U_0}"']
    &UB \otimes \lhom{UA}{\un} \ar[r,"c_{UA,UB}"']
    &\lhom{UA}{UB}.
\end{tikzcd}
$$

\begin{thm}\label{thm-char-quasi} Let $A$ be a quasi-bialgebra over a field $\kk$.
Then $A$ is a quasi-Hopf algebra if and only if the forgetful functor $U\colon A\ti\Mod \to \Vect_\kk$ admits a slack left closed structure $\mathfrak{E}$ such that $(U,\mathfrak{E})$ preserves comparison morphisms.
\end{thm}

\begin{proof} Consider the bijection between slack left Hopf structures on $A$ and slack left closed structures on $U$ given by Theorem~\ref{thm-ll-mon-mod}.
The slack left closed structure on $U$ corresponding with a slack left Hopf structure $v$ on $A$ preserves comparison morphisms if and only if for any pair of $A$-modules $V, W$, $c_{V,W}$ is the canonical map $W \otimes V^*\to \Hom_\kk(V,W)$, $x \otimes f \mapsto x \circ f$, which  means $(x\circ f) \cdot \SL(v) = (x \circ f)$ by Lemma~\ref{lemma-c-d}. Over a field this condition is equivalent to  $\SL(v) = 1 \otimes 1$, that is, $v$ is a left Hopf structure according to Theorem~\ref{thm-crit-quanti}.
\end{proof}
\subsection{Notes on the profinite case}
There is a categorical characterisation of quasi-Hopf algebras in terms of left duals in the profinite case.

A functor $U\colon \D \to \C$ between left autonomous categories is \emph{left autonomous} when it is endowed with a natural isomorphisms $U(\ldual{A}) \cong\ldual{U(A)}$ for $A$ in $\D$.

The categories of finite dimensional $A$-modules and finite dimensional vector spaces will be denoted by $A\ti\mod$ and $\vect_\kk$.

\begin{thm}\label{thm:quasi-auton}
Let $A$ be a quasi-bialgebra over a field $\kk$.
\begin{enumerate}
\item
If $A$ is a quasi-Hopf algebra, then the category $A\ti\mod$ of finite dimensional $A$-module is left autonomous, and the forgetful functor $U_0: A\ti\mod \to \vect_\kk$ is left autonomous;
\item If $A$ is profinite,
$A\ti\mod$ is left autonomous and $U_0\colon A\ti\mod \to \vect_\kk$ is left autonomous, then $A$ is a quasi-Hopf algebra.
\end{enumerate}
\end{thm}

\begin{proof} Assertion (1) goes back to Drinfeld~\cite{QHA}. The proof of Assertion (2) is similar to that of the analoguous statement for coquasi-bialgebras (see Remark~\ref{rk:coquasi}), but not clearly deductible from the latter. Sending an element of $A$ to its action on finite dimensional $A$-module defines an canonical algebra morphism $A \to \End(U_0)$, which is an isomorphism for  $A$ profinite.

Now assume $U_0$ is left autonomous. We may assume by transport of structure that the autonomous structure is an identity, that is, for each finite dimensional $A$-module $V$, $V^*$ is endowed with a action of $A$ which makes it a left dual of $V$, with the appropriate naturality conditions. This action of $A$ on $V^*$ defines an antimorphism of algebras $A \to \End(U_0)$, hence an antimorphism of algebras $\Es : A \to A$ such that the action of $A$ on $V^*$ is given by $a.f(x) = f(\Es(a)x)$.

The evaluation $e_V \colon \ldual{V} \otimes V \to \kk$ and the coevaluation $h_V \colon \kk \to \ldual{V} \otimes V$ are of the form $e_V(f \otimes x) = f(\alf_V x)$ and $h_V(1)= \Sigma_i \bee_V e_i \otimes e^i$, where $\alf_V, \bee_V : V \to V$ are natural in $V$, so by surjectivity of $A \to \End(U_0)$ there exist $\alf, \bee \in A$ such that $\alf_V$ and $\bee_V$ are the actions of $\alf$ and $\bee$ respectively.

The fact that $(\Es,\alf,\bee)$ satisfies Axioms (1)---(4) of a quasi-antipode can be derived by writing down, for $V$ in $A\ti\mod$, the relations expressing the $A$-linearity of  $e_V$ and $h_V$ and the duality relations $$(V \otimes e_V)a_{V,\ldual{V},V}(h_V \otimes V) = \id_V \blabla{ and}
(e_V \otimes \ldual{V}) a^{-1}_{\ldual{V},V,\ldual{V}} (\ldual{V} \otimes h_V) = \id_{\ldual{V}},$$ in that order, and using injectivity of $A\to  \End(U_0)$.
\end{proof}

\begin{rk}\label{rk:coquasi}
Assertion (2) is mentioned, but not clearly stated, by Majid in \cite{zbMATH00821242}, Section~{9.4}, under the hypothesis that $A$ satisfies a certain representability condition (9.45) which holds in the finite dimensional case. 
\end{rk}

Theorem~\ref{thm-char-quasi} can be deduced from
Theorem~\ref{thm:quasi-auton} in the profinite case as follows. It is well-known that  a monoidal left autonomous category $\C$ is left closed, with internal Homs $\lhom{A}{B} = B \otimes \ldual{A}$, and a monoidal left closed category $\C$ is
left autonomous if and only if for any object $V$ the comparison morphism $c_{V,V}$ is an isomorphism, and in that case $\ldual{V} = \lhom{V}{\un}$. If $v$ is a slack left Hopf structure on $A$ which preserves comparison morphisms  then $c_{V,V}$ is an isomorphism, so that $A\ti\mod$ is left autonomous, and besides, the slack left closed structure defined by $v$ provides a left autonomous structure on $U_0$, which implies that $A$ is quasi-Hopf by Theorem~\ref{thm:quasi-auton}.


\subsection{The set of slack left Hopf structures}

In a quasi-Hopf algebra $A$ not all slack left Hopf structures are left Hopf structures.  Indeed, the group $\G^l(A\otimes ?) = (A^e)^\inv$ acts freely and transitively on the right on the set of slack left Hopf structures.
Consider $A^\inv$ as a subgroup of $(A^e)^\inv$ via $t \mapsto 1 \otimes t$.  Drinfeld has shown that if $A$ is quasi-Hopf, quasi-antipodes form a free orbit under the right action (`twist') of $A^\inv$
given by  $(\Es,\alf,\bee).t = (t^{-1}\Es t, t^{-1} \alf, \bee t)$; consequently left Hopf structures form a free orbit under the right action of $A^\inv$, that is, $v \lhd t = v^{(1)} \otimes v^{(2)}t$.

On the other hand,  $(A^e)^\inv$ is the semidirect product $A^\inv \ltimes\G_0(A)$, where $$\G_0(A) = \{\gamma \in (A^e)^\inv\mid (\eps \otimes A)(\gamma) = 1\}.$$ The following theorem summarises this situation.


\begin{thm}\label{thm:quanti} Let $A$ be a quasi-Hopf algebra. The set of slack left Hopf structures of $A$ forms a torsor under the right action of $(A^e)^\inv$, inside of which the subset of left Hopf structures forms an $A^\inv$-orbit. Moreover the action of $A^\inv$ on left Hopf structures corresponds with the twist action of $A^\inv$ on the set of quasi-antipodes.

As a result for any slack left Hopf $v \in A \otimes A$, there exists a unique left Hopf structure $v_0$ and a unique element $\gamma \in \G_0(A)$ such that $v = v_0\lhd \gamma$.

Moreover in this expression of $v$, $\gamma$ is the slackness $\SL(v)$ of $v$ and $v_0$ is the left Hopf structure associated to the quasi-antipode $(\sigma, \alf, \bee')$, where $\bee' = \overline{\gamma}^{(1)}\bee\, \overline{\gamma}^{(2)}$, $\overline\gamma$ being the inverse of $\gamma$ in $A^e$.
\end{thm}

\begin{proof} The first two paragraphs result immediately from the previous remarks and the general fact that slack left Hopf structures form an orbit under the free action of $(A^e)^\inv = A^\inv \ltimes\G_0(A)$ and left Hopf structures form an orbit under  $A^\inv \subset (A^e)^\inv$.

As to the `moreover' part:
let $v =v_0 \lhd \gamma$, $v_0$ being a left Hopf structure (\emph{i.e.} $\SL(v_0) = 1 \otimes 1$) and $\gamma$ an element of $\G_0(A)$ (\emph{i.e.} $(\eps\otimes A) \gamma = 1$), with inverse $\overline{\gamma}$.
By Lemma~\ref{lem-v-change}, $\SL(v) = \SL(v_0) \cdot \gamma = \gamma$, and by Theorem~\ref{thm-quanti-from-slh}, $v_0$ is the left Hopf structure associated with the quasi-antipode $(\sigma_0,\alf_0,\bee_0)$, denoting with a subscript $0$ the data related to $v_0$, which we now compute.  By Lemma~\ref{lemma-c-d} we have $\nabla = \overline \gamma \cdot \nabla_0 \cdot \gamma$ and so $\sigma = \sigma_0$.  We also have $H^v = H^{v_0}(\gamma \,\cdot \,?)$ and so $\bee_0 = (\eps \otimes A)H^v(\overline\gamma) = \overline\gamma^{(1)} \bee_0 \overline\gamma^{(2)}$, and $\alf_0 = (\eps\otimes A)H^{v_0 -1}(1 \otimes 1) = (\eps \otimes A)(\gamma \cdot w_0) = ((\eps \otimes A)\gamma) \alf_0 = \alf_0$.
%
%
%
%
\end{proof}

\begin{rk}
The assumption we made in this Section that the base ring $k$ is a field, which is used for the first time in the proof of Lemma~\ref{lemma:sigma-delta},  may be lifted provided $A$ is assumed to be projective as a $k$-module.
\end{rk}

\section{Conclusion}\label{sec:concl}

We wish to discuss briefly  two questions which remain open. The first is, does slack Hopf imply Hopf?
The question, so stated, is a bit vague, since the prefix Hopf doesn't have the same meaning for, say, a bialgebra and a quasi-bialgebra, so let's restrict ourselves to comonoidal monads on monoidal categories, or bialgebras over a field.

For a comonoidal monad $T$ to be left Hopf seems a much stronger assumption than to be slack left Hopf since in the first case, a certain natural transformation $F^T(X \otimes U^T A) \to F^T(X) \otimes A$ must be an isomorphism, whereas in the second it is only required that there exist an arbitrary natural isomorphism between those functors. However at this point we don't dispose of any example of a slack left Hopf comonoidal monad which isn't left Hopf.

In the following cases we know that slack Hopf comonoidal monads (or monoidal comonads) are Hopf:
\begin{itemize}
    \item on a cartesian category (Theorem~\ref{th-hopf-cart});
    \item for a commutative algebra (Lemma~\ref{lemma:commut});
    \item for the monoidal comonad associated to a small category (Theorem~\ref{thm-cat-comonad}).
\end{itemize}

These three cases have in common that the methods we used to prove `slack Hopf implies Hopf' are very specific and apparently not significantly generalisable. Indeed, in the first case, the proof rests on the fact that in a cartesian category a slack Hopf structure $\beta$ is necessarily split, \emph{i.e.} of the form $\beta = \gamma \times \delta$, with $\gamma, \delta \in \Nat(\id_\C,T)$. In the second case, it uses the fact that for a commutative bialgebra the expression of the morphism $H^v$ can be factorised in a manner which is impossible for a general bialgebra (see Remark~\ref{rk:near-com}). In the third example, it comes from the fact that owing to the rigidity of the situation there is only one possible candidate for being a slack Hopf structure (Lemma~\ref{lemma:rigid}).

Note that it is not difficult to construct a slack Hopf colax magma monad for which the fusion operator $H^{\eta \otimes \eta}$ is not an isomorphism, such as that of Example~\ref{exa:mat} or certain quasi-Hopf algebras, but those examples are not comonoidal.

The second question is: what would be a quasi-Hopf monad?
We've seen that quasi-bialgebras are slack Hopf, but we don't know yet whether the converse is true. We have however a characterization of quasi-bialgebras in terms of slack Hopf structures, namely a quasi-Hopf algebra is a quasi-bialgebra having a slack Hopf structure $v$ whose slackness $\SL(v)$ is invertible, or, equivalently, is trivial (Corollary~\ref{cor-quanti}).  This algebraic criterion can be reformulated in categorical terms: a quasi-Hopf algebra is a quasi-bialgebra $A$ whose forgetful functor $A\ti\Mod \to \Vect_\kk$ has a slack left closed structure which preserves the comparison morphisms (Theorem~\ref{thm-char-quasi}).

The existence of this purely categorical criterion suggests that the notion of a quasi-Hopf algebra may be generalised to a notion of a quasi-Hopf monad, a possibility  we plan to explore in a subsequent work.


\begin{thebibliography}{10}

\bibitem{NonAss}
G.~E. Barnes, A.~Schenkel, and R.~J. Szabo.
\newblock Nonassociative geometry in quasi-{Hopf} representation categories.
  {I}: {Bimodules} and their internal homomorphisms.
\newblock {\em J. Geom. Phys.}, 89:111--152, 2015.

\bibitem{zbMATH02222246}
M.~Barr and C.~Wells.
\newblock Toposes, triples and theories.
\newblock {\em Repr. Theory Appl. Categ.}, 2005(12):1--288, 2005.

\bibitem{zbMATH06450114}
G.~B{\"o}hm, Y.~Chen, and L.~Zhang.
\newblock On {Hopf} monoids in duoidal categories.
\newblock {\em J. Algebra}, 394:139--172, 2013.

\bibitem{MR2793022}
A.~Brugui{\`e}res, S.~Lack, and A.~Virelizier.
\newblock Hopf monads on monoidal categories.
\newblock {\em Adv. Math.}, 227(2):745--800, 2011.

\bibitem{BruVir}
A.~Brugui{\`e}res and A.~Virelizier.
\newblock Hopf monads.
\newblock {\em Advances in Mathematics}, 215(2):679--733, 2007.

\bibitem{BV3}
A.~Brugui{\`e}res and A.~Virelizier.
\newblock Quantum double of hopf monads and categorical centers.
\newblock {\em Transactions of the Americal Mathematical Society},
  364(3):1225--1279, 2012.

\bibitem{Day-Street:quantumcat}
B.~Day and R.~Street.
\newblock Quantum categories, star autonomy, and quantum groupoids.
\newblock In {\em Galois theory, Hopf algebras, and semiabelian categories},
  volume~43 of {\em Fields Inst. Commun.}, pages 187--225. Amer. Math. Soc.,
  Providence, RI, 2004.

\bibitem{QHA}
V.~G. Drinfeld.
\newblock Quasi-hopf algebras.
\newblock {\em Leningrad Math. J.}, 1(6):1419–--1457, 1990.

\bibitem{kelly_doctrinal_1974}
G.~M. Kelly.
\newblock Doctrinal adjunction.
\newblock In {\em Category {Seminar} ({Proc}. {Sem}., {Sydney}, 1972/1973)},
  Lecture Notes in Math., Volume 420, pages 257--280. Springer, Berlin, 1974.

\bibitem{zbMATH02172008}
G.~M. Kelly.
\newblock Basic concepts of enriched category theory.
\newblock {\em Repr. Theory Appl. Categ.}, 2005(10):1--136, 2005.

\bibitem{zbMATH03223713}
H.~Kleisli.
\newblock Every standard construction is induced by a pair of adjoint functors.
\newblock {\em Proc. Am. Math. Soc.}, 16:544--546, 1965.

\bibitem{CWM}
S.~Mac~Lane.
\newblock {\em Categories for the working mathematician}.
\newblock Springer-Verlag, New York, second edition, 1998.

\bibitem{zbMATH00821242}
S.~Majid.
\newblock {\em Foundations of quantum group theory}.
\newblock Cambridge: Cambridge Univ. Press, 1995.

\bibitem{zbMATH01863401}
P.~McCrudden.
\newblock Opmonoidal monads.
\newblock {\em Theory Appl. Categ.}, 10:469--485, 2002.

\bibitem{Moer}
I.~Moerdijk.
\newblock Monads on tensor categories.
\newblock {\em J. Pure Appl. Algebra}, 168(2--3):189–--208, 2002.

\bibitem{Schau:CFQHA}
P.~Schauenburg.
\newblock Two characterizations of finite quasi-{Hopf} algebras.
\newblock {\em J. Algebra}, 273(2):538--550, 2004.

\end{thebibliography}
\end{document}